\documentclass[11pt]{article}
\usepackage{amsmath,amsthm,amssymb,xypic}

\theoremstyle{plain}  

\newtheorem{thm}{Theorem}[section]
\newtheorem{prop}[thm]{Proposition}
\newtheorem{lem}[thm]{Lemma}

\theoremstyle{definition}

\newtheorem{prop-defn}[thm]{Proposition--Definition}
\newtheorem{defn}[thm]{Definition}
\newtheorem{notn}[thm]{Notation}

\theoremstyle{remark}

\newtheorem{rem}[thm]{Remark}
\newtheorem{ex}[thm]{Example}

\DeclareMathOperator{\Aut}{Aut}
\DeclareMathOperator{\Cl}{Cl}

\DeclareMathOperator{\Spec}{Spec}

\DeclareMathOperator{\Hom}{Hom}
\DeclareMathOperator{\cHom}{\mathcal{H}\mathnormal{om}}
\DeclareMathOperator{\Ext}{Ext}
\DeclareMathOperator{\cExt}{\mathcal{E}\mathnormal{xt}}

\DeclareMathOperator{\Diff}{Diff}
\DeclareMathOperator{\Pic}{Pic}

\DeclareMathOperator{\Supp}{Supp}
\DeclareMathOperator{\Exc}{Ex}

\DeclareMathOperator{\wt}{wt}

\DeclareMathOperator{\nd}{\not \: \mid}

\newcommand{\QED}{\ifhmode\unskip\nobreak\fi\quad {\rm Q.E.D.}} 

\newcommand{\bA}{\mathbb A}
\newcommand{\bC}{\mathbb C}
\newcommand{\bF}{\mathbb F}

\newcommand{\bN}{\mathbb N}
\newcommand{\bP}{\mathbb P}
\newcommand{\bQ}{\mathbb Q}
\newcommand{\bR}{\mathbb R}

\newcommand{\bZ}{\mathbb Z}

\newcommand{\cB}{\mathcal B}

\newcommand{\cX}{\mathcal X}
\newcommand{\cY}{\mathcal Y}
\newcommand{\sD}{\mathcal D}

\newcommand{\cF}{\mathcal F}
\newcommand{\cG}{\mathcal G}

\newcommand{\cI}{\mathcal I}
\newcommand{\cO}{\mathcal O}
\newcommand{\sL}{\mathcal L}
\newcommand{\cM}{\mathcal M}

\newcommand{\cT}{\mathcal T}
\newcommand{\cU}{\mathcal U}

\newcommand{\cZ}{\mathcal Z}

\newcommand{\stX}{\mathfrak X}
\newcommand{\stD}{\mathfrak D}

\newcommand{\map}{\rightarrow}
\newcommand{\da}{\downarrow}
\newcommand{\inj}{\hookrightarrow}

\title{Compact moduli of plane curves}
\author{Paul Hacking}
\date{July 12, 2003}
\begin{document}
\maketitle
\pagestyle{plain}

\begin{abstract}
We construct a compactification $\cM_d$ of the moduli space of plane curves of degree $d$.
We regard a plane curve $C \subset \bP^2$ as a surface-divisor pair $(\bP^2,C)$ and 
define $\cM_d$ as a moduli space of pairs $(X,D)$ where $X$ is a degeneration of the plane.
We show that, if $d$ is not divisible by 3, the stack $\cM_d$ is smooth and 
the degenerate surfaces $X$ can be described explicitly.

\end{abstract}

\section{Introduction}
Let $V_d$ be the moduli space of smooth plane curves of degree $d \ge 3$.
Then $V_d$ is the quotient $U_d/ \Aut(\bP^2)$ where $U_d$ is the
open locus of smooth curves in the Hilbert scheme $H_d$ of plane curves
of degree $d$. These moduli spaces are fundamental objects in algebraic
geometry. Geometric invariant theory provides a compactification $\bar{V_d}$
of $V_d$. However, $\bar{V_d}$ is rather unsatisfactory for several reasons.
First, $\bar{V_d}$ is not a moduli space itself --- some points of the boundary
correspond to several isomorphism classes of plane curves. Second,
$\bar{V_d}$ has fairly complicated singularities at the boundary.
In particular, these rule out the possibility of performing 
intersection theory on $\bar{V_d}$
to obtain enumerative results. Finally, the boundary is difficult to describe
explicitly --- there is a stratification  given by the type of 
singularities on the
degenerate curve, but this can only be computed for small degrees.

In this paper we describe an alternative compactification $\cM_d$ of $V_d$.
The space $\cM_d$ is a moduli space of \emph{stable pairs}.
A stable pair is a surface-divisor pair $(X,D)$  
which is a degeneration of the plane together with a curve and 
satisfies certain additional properties. 
Morally speaking, the pair $(X,D)$ should be identified with the curve $D$;
the existence of an embedding $D \inj X$ gives some structural information
about $D$, e.g., the existence of a Brill--Noether special 
linear system on $D$. There is a stratification of $\cM_d$ given by the 
isomorphism type of the surface $X$.
If $d$ is not divisible by 3 then we can explicitly describe the surfaces $X$ 
which occur and so determine this stratification. 
Moreover, in this case, the space $\cM_d$ is smooth (as a stack) and,
writing $\cM^0_d$ for the open stratum corresponding to the plane, the 
boundary $\cM_d \backslash \cM^0_d$ is a normal crossing divisor.
 
We pause to describe the simplest example, namely the case $d=4$.
The surfaces $X$ occuring are the plane, 
the cone over the rational normal curve of degree $4$ and the 
non-normal surface obtained by glueing two 
quadric cones along a ruling so that the vertices coincide.
In the language of weighted projective spaces, the latter two surfaces are 
$\bP(1,1,4)$ and $\bP(1,1,2) \cup \bP(1,1,2)$ respectively.
The curves lying on $\bP(1,1,4)$ are hyperelliptic ---
we obtain a 2-to-1 map to $\bP^1$ by projecting away from the vertex. 
The curves lying on $\bP(1,1,2) \cup \bP(1,1,2)$ are 
`degenerate hyperelliptic' 
--- we obtain a 2-to-1 map to $\bP^1 \cup \bP^1$ by projecting away 
from the common vertex of the component surfaces; 
these curves have two components of genus 1 meeting in two nodes. 
The stratification of $\cM_4$ is as follows: we have
$\cM_4 = Z_0 \cup Z_1 \cup Z_2$ where  $Z_0$, $Z_1$ and $Z_2$ denote
the strata corresponding to $\bP^2$, $\bP(1,1,4)$ and 
$\bP(1,1,2) \cup \bP(1,1,2)$.
The stratum $Z_0$ is open, 
$Z_1$ is a locally closed locus of codimension 1, $Z_2$ is closed of 
codimension 2 and the closure of $Z_1$ is $Z_1 \cup Z_2$. 
The degree $4$ case was originally treated by Hassett \cite{Has}, who worked
with a different class of pairs $(X,D)$. 
Roughly, we allow worse singularities on $D$ in order to gain greater control
of the surface $X$. The amazing thing is that, with our definition of 
stable pair, many of the features of the degree 4 case persist for all degrees
which are not divisible by 3 (e.g. $\cM_d$ is smooth and each degenerate 
surface $X$ has at most two components).

We describe stable pairs in more detail. 
If $(X,D)$ is a stable pair then the surface $X$ has semi log canonical 
singularities (Definition~\ref{defn-slc}) and the $\bQ$-Cartier divisor $-K_X$
is ample. The divisor $D$ lies in the linear system $|\frac{-d}{3}K_X|$ and has
mild singularities. More carefully, the singularities of $D$ 
which are permitted are precisely those such that the log canonical threshold 
of the pair $(X,D)$ is strictly larger than $\frac{3}{d}$. For example, if
$d=4$, the singularities of $D$ are either nodes or cusps.

There is a coarse classification of the surfaces $X$ into types A, B, C and D.
Type A are the normal surfaces. Type B have two normal components meeting in a
 smooth rational curve. Types C and D have several components  forming an
`umbrella' or a `fan'  respectively. If the degree $d$ is not divisible by 3 
then only types A and B occur; in particular, $X$ has at most 2 components.
Moreover, the only singularities of $X$ are quotients of smooth or normal 
crossing points.

We give an explicit description of the surfaces $X$ of type A. 
If $X$ is log terminal then $X$ is obtained as a deformation of
a weighted projective space $\bP(a^2,b^2,c^2)$ where $(a,b,c)$ is a solution 
of the Markov equation
$$a^2+b^2+c^2=3abc.$$ 
This is a refinement of a result of Manetti \cite{Ma}, so we call such 
surfaces \emph{Manetti surfaces}. If $X$ is not log terminal
then $X$ is an elliptic cone of degree 9. 
We also present a finer classification of the surfaces of type B.

We give a map of the paper. In Section~\ref{stablepairs} we define stable pairs and prove a completeness property, namely, 
that a family of smooth plane curves over a punctured curve can be completed to a family of stable pairs in a canonical way.
In Section~\ref{QG} we develop a theory of $\bQ$-Gorenstein deformations for semi log canonical surfaces which we use
to construct the moduli space of stable pairs $\cM_d$ and study its infinitesimal properties. We construct the space $\cM_d$ 
in Section~\ref{moduli}. We also provide an effective bound on the index of a surface occurring in a stable pair in terms of the degree.
In Section~\ref{coarse} we give the coarse classification of the degenerate surfaces $X$. 
In Section~\ref{prelim} we collect some restrictions on the singularities of $X$ and the
Picard numbers of the components implied by the existence of a smoothing of to $\bP^2$.
Section~\ref{simp} provides the simplifications in the case $3 \nd d$ stated above and Sections~\ref{normal} and \ref{typeB}
give the classification of the type A and  B surfaces respectively. In Section~\ref{GIT} we explain the relation between our notion of stability and 
GIT stability for a plane curve.  
Finally, in Section~\ref{Examples} we give the complete classification of stable pairs of degrees $4$ and $5$.

This paper is based on my PhD thesis \cite{Hac1} and subsequent recent work. 
I would like to thank my supervisor, Alessio Corti, for constant guidance, 
encouragement and friendship throughout the course of my PhD. I am also grateful to Brendan Hassett, S\'{a}ndor Kov\'{a}cs,
Miles Reid, and Nick Shepherd-Barron for various helpful discussions.
The final version of this article will be published in the Duke Mathematical Journal, published by
Duke University Press.

\section{Stable pairs} \label{stablepairs} 

We define the notion of a stable pair and show that, possibly after base change, every family of smooth plane curves over a punctured curve
can be completed to a family of stable pairs in a unique way.
Equivalently, the moduli space of stable pairs is separated and proper.
As a preliminary step, we define semistable pairs and show that every such family can be completed to a family of semistable pairs,
although the completion is not necessarily uniquely determined.

We use the semistable minimal model program, which is explained
in \cite{KM}, Chapter~7. Our construction is a refinement of the usual 
construction of compact moduli of pairs (\cite{KSB},\cite{Al}) 
applied to the case of pairs consisting of the plane together with a curve of 
degree $d$.
The standard construction produces a moduli space 
$\cM^{\alpha}_d$ of pairs
$(X,D)$ such that $(X, \alpha D)$ is semi log canonical (see Definition~\ref{defn-slc}) and $K_X+\alpha D$  
is ample for some fixed $\alpha \in \bQ$; here we require 
$\alpha>\frac{3}{d}$ in order that $K_{\bP^2}+\alpha D$ is ample for $D$ 
a plane curve of degree $d$.
However, there are technical problems in the construction of this 
moduli space, in particular, 
the correct definition of a family $(\cX,\sD)/S$ of such pairs is 
unclear. The main problem is that we cannot insist that both the relative 
divisors $K_{\cX}$ and $\sD$ are $\bQ$-Cartier. 
This complicates the deformation theory and thus renders an infinitesimal 
study of the moduli space intractable.
We instead construct a moduli space $\cM_d$ of \emph{stable pairs}.
A stable pair is a pair $(X,D)$ such that 
$(X,(\frac{3}{d}+ \epsilon)D)$ is semi log canonical 
and $K_X+(\frac{3}{d}+\epsilon)D$ is ample for all $0 < \epsilon \ll 1$.
It is not immediately clear that stable pairs are bounded; however, once this
is established, we deduce that $\epsilon$ may be chosen uniformly. That is,
there exists $\epsilon_0 > 0$ such that for every stable pair $(X,D)$, 
the pair $(X,(\frac{3}{d}+ \epsilon)D)$ is semi log canonical and 
$K_X+(\frac{3}{d}+\epsilon)D$ is ample for any $0 < \epsilon \le \epsilon_0$.
Thus $\cM^{\alpha}_d$ coincides with $\cM_d$ for 
$\frac{3}{d} < \alpha \le \frac{3}{d} + \epsilon_0$, or, more coarsely,
$\cM_d$ is the limit of $\cM^{\alpha}_d$ as $\alpha \searrow \frac{3}{d}$.
This was the original motivation for the definition of a stable pair.
The space $\cM_d$ is much easier to understand than the space $\cM^{\alpha}_d$
for arbitrary $\alpha$. Hence, in what follows, we construct $\cM_d$ directly.

\begin{notn}
We always work over $\bC$.
We write $0 \in T$ for the germ of a smooth curve.
We use script letters to denote flat families over $T$
and regular letters for the special fibre, e.g.,
\begin{eqnarray*}
\begin{array}{ccc}
X & \subset & \cX \\
\da &       & \da \\
0   & \in   & T
\end{array}
\end{eqnarray*}
\end{notn} 

We recall the definition of semi log canonical singularities of 
surface-divisor pairs (\cite{KSB}, \cite{Al}).
These are the singularities we must allow to compactify moduli of pairs.  

\begin{defn} \label{defn-slc}
Let $X$ be a surface and $D$ an effective $\bQ$-divisor on $X$.
The pair $(X,D)$ is semi log canonical (respectively semi log terminal) if
\begin{enumerate}
\item The surface $X$ is Cohen-Macaulay and has only normal crossing singularities in codimension $1$.
\item Let $K_X$ denote the Weil divisor class on $X$ corresponding to the dualising sheaf $\omega_X$.
Then the divisor $K_X+D$ is $\bQ$-Cartier.
\item Let $\nu \colon X^{\nu} \map X$ be the normalisation of $X$. Let $\Delta$ denote the double curve 
of $X$ and write $D^{\nu}$ and $\Delta^{\nu}$ for the inverse images of $D$ and $\Delta$ on $X^{\nu}$.
Then the pair $(X^{\nu},\Delta^{\nu}+D^{\nu})$ is log canonical (respectively log terminal).
\end{enumerate}
We use the abbreviations slc and slt for semi log canonical and semi log terminal.
\end{defn}

\begin{rem}

\begin{enumerate}
\item The dualising sheaf $\omega_X$ satisfies Serre's condition $S_2$. It is also invertible in 
codimension $1$ by (1). Hence it corresponds to a Weil divisor class $K_X$ as stated. If $X$ is normal 
this is of course the usual canonical divisor class.
\item If $(X,D)$ is slc then no component of $D$ is contained in the double curve $\Delta$ by (3).
\item Note that $K_{X^{\nu}}+\Delta^{\nu}+D^{\nu}=\nu^{\star}(K_X+D)$.
\end{enumerate}
\end{rem}

\begin{defn}\label{defn-semistable}
Let $X$ be a surface and $D$ an effective $\bQ$-Cartier divisor on $X$. Let $d \in \bN$, $d \ge 3$.
The pair $(X,D)$ is a \emph{semistable pair} of degree $d$ if
\begin{enumerate}
\item The surface $X$ is normal and log terminal.
\item The pair $(X,\frac{3}{d}D)$ is log canonical.
\item The divisor $dK_X+3D$ is linearly equivalent to zero.
\item There is a deformation $(\cX,\sD)/T$ of the pair $(X,D)$ over the germ of a curve such that
the general fibre $\cX_t$ of $\cX/T$ is isomorphic to $\bP^2$ and the divisors $K_{\cX}$ and $\sD$ 
are $\bQ$-Cartier.
\end{enumerate}
\end{defn}

\begin{rem}
There is a very concrete classification of the surfaces $X$ appearing here (Theorem~\ref{thm-manetti}).
\end{rem}

\begin{thm}\label{thm-semistable}
Let $0 \in T$ be a germ of a curve and write $T^{\times}=T-0$.
Let $\sD^{\times} \subset \bP^2 \times T^{\times}$ be a family of smooth plane curves over $T^{\times}$ of degree $d \ge 3$.
Then there exists a finite surjective base change $T' \map T$ and a family $(\cX,\sD)/T'$ of semistable pairs
extending the pullback of the family $(\bP^2 \times {T^{\times}}, {\sD^{\times}})/ T^{\times}$ 
such that the divisors $K_{\cX}$ and $\sD$ are $\bQ$-Cartier.
\end{thm}
\begin{proof}
First complete $(\bP^2 \times T^{\times}, \sD^{\times})$ to a flat family $(\bP^2 \times T, \sD)$ over $T$.
After a base change (which we will suppress in our notation) there is a semistable log resolution 
$$\pi \colon (\tilde{\cX},\tilde{\sD}) \map (\bP^2 \times T, \sD)/T$$
which is an isomorphism over $T^{\times}$. We proceed as follows: 
\begin{enumerate}
\item Run a $K_{\tilde{\cX}}+\frac{3}{d}\tilde{\sD}$ MMP over $T$. 
Let $(\cX_1,\sD_1)/T$ denote the end product. Then $K_{\cX_1}+\frac{3}{d}\sD_1$ is relatively nef and 
vanishes on $\cX_1^{\times}=\bP^2 \times T^{\times}$; 
it follows that $dK_{\cX_1}+3\sD_1 \sim 0$ by Lemma~\ref{kulikov}.
\item Run a $K_{\cX_1}$ MMP over $T$. The end product $(\cX,\sD)/T$ is the required completion of 
 $(\bP^2 \times T^{\times}, \sD^{\times})$.
\end{enumerate}
We verify the required properties of $(\cX,\sD)/T$.
We refer to \cite{KM} Chapter~7 for background on the semistable minimal model program.
The family $\cX/T$ is a Mori fibre space since it is the end product of a MMP and 
the general fibre is a del Pezzo surface, namely $\bP^2$. 
Regarding the singularities of $\cX/T$, we know that the pair $(\cX,X)$ is dlt and $\cX$
is $\bQ$-factorial. It follows that $X$ is irreducible using $\rho(\cX/T)=1$ and $\bQ$-factoriality.
Then $X$ is normal and log terminal by the dlt property.

The pair $(\cX_1,X_1+\frac{3}{d}\sD_1)$ is dlt; since 
$dK_{\cX_1}+3\sD_1 \sim 0$ it follows that $dK_{\cX}+3\sD \sim 0$ and
$(\cX,X+\frac{3}{d}\sD)$ is log canonical. Thus $(X,\frac{3}{d}D)$
is log canonical and $dK_X+3D \sim 0$ by adjunction.
\end{proof} 

\begin{lem} \label{kulikov}
Let $\cX/(0 \in T)$ be a flat family of projective slc surfaces over the germ of a curve 
such that the general fibre is normal. 
Let $\cX^{\times}/T^{\times}$ denote the restriction 
of the family to the punctured curve $T^{\times}= T \backslash \{ 0 \}$.
Let $\cB$ be a $\bQ$-Cartier divisor on $\cX$ such that $\cB$ is relatively nef
and $\cB|_{\cX^{\times}} \sim 0$. Then $\cB \sim 0$.
\end{lem}
\begin{proof}
Let $X_1, \ldots , X_n$ denote the irreducible components of $X$, so $X = \sum X_i$
as divisors on $\cX$.
We have an exact sequence
$$0 \map \bZ X \map \oplus \bZ X_i \map \Cl (\cX) \map \Cl (\cX^{\times}) \map 0.$$
Hence, since $\cB|_{\cX^{\times}} \sim 0$, we may write 
$\cB \sim \sum a_iX_i$, where $a_i \le 0$ for all $i$ and we have equality for some $i$.
If $a_j=0$, then $\cB|_{X_j}=\sum_{i \neq j} a_iX_i|_{X_j} \le 0$. But $\cB|_{X_j}$ is nef,
hence $\cB|_{X_j}=0$, i.e., $a_i=0$ for each $i$ such that $X_i$ and $X_j$ meet in a curve.
It follows by induction that $a_i=0$ for all $i$, i.e., $\cB  \sim 0$.
\end{proof}

\begin{defn}\label{defn-stable}
Let $X$ be a surface and $D$ an effective $\bQ$-Cartier divisor on $X$. Let $d \in \bN$, $d \ge 4$.
The pair $(X,D)$ is a \emph{stable pair} of degree $d$ if
\begin{enumerate}
\item The pair $(X,(\frac{3}{d}+\epsilon)D)$ is slc and the divisor $K_X+(\frac{3}{d}+\epsilon)D$ is ample 
for some $\epsilon >0$.
\item(=\ref{defn-semistable}(3)) The divisor $dK_X+3D$ is linearly equivalent to zero.
\item(=\ref{defn-semistable}(4)) There is a deformation $(\cX,\sD)/T$ of the pair $(X,D)$ over the germ of a curve such that
the general fibre $\cX_t$ of $\cX/T$ is isomorphic to $\bP^2$ and the divisors $K_{\cX}$ and $\sD$ 
are $\bQ$-Cartier.
\end{enumerate}
\end{defn}

\begin{rem} Conditions (1) and (2) may be replaced by the following 
(cf. our motivating remarks in the introduction of this section):
\begin{enumerate} 
\item[($1'$)] The pair $(X,(\frac{3}{d}+\epsilon)D)$ is slc and the divisor $K_X+(\frac{3}{d}+\epsilon)D$ is ample 
for all $0 < \epsilon \ll 1$.
\end{enumerate}
Clearly (1) and (2) imply ($1'$) and ($1'$) implies (1); it remains to show that ($1'$) (together with (3)) implies (2).
If $(X,D)$ satisfies ($1'$) then, 
since $K_X+(\frac{3}{d} +\epsilon)D$ is ample for all $0 < \epsilon \ll 1$, the limit $K_X+\frac{3}{d}D$ is nef.
Suppose $(\cX,\sD)/T$ is a smoothing of $(X,D)$ as in (3). The divisor $dK_{\cX}+3\sD$ is relatively nef
and vanishes on the general fibre, hence is linearly equivalent to zero by Lemma~\ref{classgp}(1) and 
Lemma~\ref{kulikov}. Thus
$dK_X+3D \sim 0$ by restriction, so $(X,D)$ satisfies (2) as required.
\end{rem}

\begin{rem} 
We note that, if $d$ is a multiple of 3, then $\frac{d}{3}K_X+D \sim 0$.
For, writing $(\cX,\sD)/T$ for a smoothing as above, the condition $dK_X+3D \sim 0$ implies that $dK_{\cX}+3\sD \sim 0$
and $\Cl(\cX)$ is torsion-free by Lemma~\ref{classgp}, hence $\frac{d}{3}K_{\cX}+\sD \sim 0$ and
so $\frac{d}{3}K_X+D \sim 0$ by restriction.
\end{rem} 

\begin{lem} \label{classgp}
Let $\cX/(0 \in T)$ be a flat family of surfaces over the germ of a curve 
with general fibre $\bP^2$ and reduced special fibre $X$. Then 
\begin{enumerate}
\item $\cX^{\times} \cong \bP^2 \times T^{\times}$
\item $\Cl(\cX) \cong \bZ^n$, where $n$ is the number of components of $X$.
\end{enumerate}
\end{lem}
\begin{proof}
Since the general fibre is $\bP^2$ there is no monodromy and $\cX^{\times} \cong \bP^2 \times T^{\times}$.
Hence $\Cl(\cX^{\times}) \cong \bZ$. The exact sequence
$$0 \map \bZ X \map \oplus \bZ X_i \map \Cl (\cX) \map \Cl (\cX^{\times}) \map 0$$
now gives $\Cl(\cX) \cong \bZ^n$ as claimed.
\end{proof}

\begin{thm}\label{thm-stable}
Let $\sD^{\times} \subset \bP^2 \times T^{\times}$ be a family of smooth plane curves of degree $d \ge 4$ over a punctured curve 
$T^{\times}$.
Then there exists a finite surjective base change $T' \map T$ and a family $(\cX,\sD)/T'$ of stable pairs
extending the pullback of the family $(\bP^2 \times T^{\times}, \sD^{\times})/ T^{\times}$ such that the divisors $K_{\cX}$ and $\sD$ 
are $\bQ$-Cartier. Moreover the family $(\cX,\sD)/T'$ is unique in the following sense: 
any two such families become isomorphic after a further finite surjective base change. 
\end{thm}
\begin{proof} 
Let $(\cX_1,\sD_1)/T$ be a family of semistable pairs extending the family $(\bP^2 \times T^{\times}, \sD^{\times})/ T^{\times}$
as constructed in the proof of Theorem~\ref{thm-semistable}. Then the pair $(\cX_1,X_1+\frac{3}{d}\sD_1)$ is log canonical and
the pair $(\cX_1,X_1)$ is dlt. 
There exists a partial semistable resolution (a `maximal crepant blowup' of $(\cX_1,\frac{3}{d}\sD_1)$)
$$\pi \colon (\cX_2,\sD_2) \map (\cX_1,\sD_1)/T$$ 
such that $dK_{\cX_2}+3\sD_2=\pi^{\star}(dK_{\cX_1}+3\sD_1) \sim 0$ and 
$(\cX_2,X_2+(\frac{3}{d}+\epsilon)\sD_2)$ is dlt for $0< \epsilon \ll 1$.
Let $(\cX,\sD)/T$ be the $K_{\cX_2}+(\frac{3}{d}+\epsilon)\sD_2$ canonical model. 
Then $(\cX,X+(\frac{3}{d}+\epsilon)\sD)$ is log canonical, the divisor $dK_{\cX}+3\sD \sim 0$ and $K_{\cX}+X+(\frac{3}{d}+\epsilon)\sD$
is relatively ample.
By adjunction $(X,(\frac{3}{d}+\epsilon)D)$ is slc, the divisor $dK_X+3D \sim 0$ and $K_X+(\frac{3}{d}+\epsilon)D$ is ample.
Note also that $K_{\cX}+(\frac{3}{d}+\epsilon)\sD$ is $\bQ$-Cartier by construction.
Hence $K_{\cX}$ and $\sD$ are $\bQ$-Cartier since $dK_{\cX}+3\sD \sim 0$.

To prove uniqueness, note that $(\cX,\sD)/T$ is the $K_{\tilde{\cX}}+(\frac{3}{d}+\epsilon)\tilde{\sD}$ canonical model of any semistable 
log resolution  $(\tilde{\cX},\tilde{\sD})/T$, where $\epsilon>0$ is sufficiently small.
\end{proof}

We record the following important result, which is an immediate consequence of 
conditions (1) and (2) of Definition~\ref{defn-stable}.
\begin{prop} 
Let $(X,D)$ be a stable pair. Then $X$ is an slc surface and the divisor $-K_X$ is ample.
\end{prop}

\section{$\bQ$-Gorenstein deformation theory} \label{QG}

We define the $\bQ$-Gorenstein deformations of a slc surface $X$ 
to be those locally induced by a deformation of the canonical covering of $X$. 
We then describe how to calculate the $\bQ$-Gorenstein deformations of a given surface $X$.
This theory is used in Section~\ref{moduli} to construct the moduli space $\cM_d$ of stable pairs
and in Section~\ref{simp} to prove that $\cM_d$ is smooth if $3 \nd d$.
It can also be used to construct compact moduli spaces of surfaces of general type
with a finer scheme theoretic structure than that originally defined in \cite{KSB}
and facilitates an infinitesimal study of such moduli spaces.
My presentation here is influenced by earlier work of Koll\'{a}r and Hassett \cite{Has}.

If a sheaf $\cF$ on a surface $X$ satisfies the $S_2$ condition, one can recover $\cF$ from $\cF|_U$ 
where $U \inj X$ has finite complement.
We require a relative $S_2$ condition for sheaves on families of slc surfaces 
which allows us to do this in the relative context. 
The definition and basic results are collected in Appendix~\ref{app_S_2}.

\subsection{Definition of $\bQ$-Gorenstein deformations}

Let $P \in X$ be an slc surface germ.
We define the \emph{canonical covering} $\pi \colon Z \map X$ by
$$Z = \underline{\Spec}_X(\cO_X \oplus \cO_X(K_X) \oplus \cdots \oplus \cO_X((N-1)K_X)),$$  
where $N$ is the index of $P \in X$ and the multiplication is given by fixing an isomorphism 
$\cO_X(NK_X) \stackrel{\sim}{\map} \cO_X$. This is a straightforward generalisation of the usual construction for
$X$ a normal variety such that $K_X$ is $\bQ$-Cartier (cf. \cite{YPG}).
It is characterised by the following properties:
\begin{enumerate}
\item The morphism $\pi$ is a cyclic quotient of degree $N$ which is \'{e}tale in codimension $1$.
\item The surface $Z$ is Gorenstein, i.e., it is Cohen-Macaulay and the Weil divisor $K_Z$ is Cartier.
\end{enumerate}

For $X$ an slc surface, the canonical covering at a point $P \in X$ is uniquely determined in the \'{e}tale topology. 
Hence the data of canonical coverings everywhere locally on $X$ defines a Deligne-Mumford stack $\stX$ with 
coarse moduli space $X$, the \emph{canonical covering stack} of $X$ (cf. \cite{Ka}, p. 18, Definition~6.1).

\begin{defn}
Let $P \in X$ be an slc surface germ. 
Let $N$ be the index of $X$ and $Z \map X$ the canonical covering, a $\mu_N$ quotient.
We say a deformation $\cX/(0 \in S)$ of $X$ is \emph{$\bQ$-Gorenstein} if
there is a $\mu_N$-equivariant deformation $\cZ/S$ of $Z$ whose quotient is $\cX/S$.
\end{defn}

\begin{notn}
Let $\cX/S$ be a flat family of slc surfaces. Let $i \colon \cX^0 \inj \cX/S$ be the inclusion of the Gorenstein locus
of $\cX/S$, i.e., the locus where the relative dualising sheaf $\omega_{\cX/S}$ is invertible. We 
write $\omega_{\cX/S}^{[N]}$ for the sheaf $i_{\star}\omega_{\cX^0/S}^{\otimes N}$.
\end{notn} 

We say that a family $\cX/S$ is \emph{weakly $\bQ$-Gorenstein} if the sheaf $\omega_{\cX/S}^{[N]}$ is invertible
for some $N \ge 1$ (cf. \cite{KSB}). The least such $N$ is the \emph{index} of $\cX/S$. 
If $\cX$ is normal and $S$ is smooth this is just the requirement that $K_{\cX}$ is $\bQ$-Cartier.
We show that a $\bQ$-Gorenstein family is weakly $\bQ$-Gorenstein (Lemma~\ref{lem-weakQG}).
Moreover, if the base $S$ is a curve and the general fibre of $\cX/S$ is canonical then the two conditions are equivalent
(Lemma~\ref{lem-DVR}). 

\begin{lem} \label{lem-weakQG}
Let $P \in X$ be an slc surface germ of index $N$. 
Let $\cX/(0 \in S)$ be a $\bQ$-Gorenstein deformation of $X$.
Then $\cX/S$ is weakly $\bQ$-Gorenstein of index $N$.
\end{lem}
\begin{proof}
There is a diagram
\begin{eqnarray*}
\begin{array}{ccc}
Z & \subset & \cZ \\
\da &       & \da \\
X & \subset & \cX \\
\da &       & \da \\
0   & \in   & S
\end{array}
\end{eqnarray*}
where $Z$ is the canonical cover of $P \in X$ and $\cZ/S$ is a $\mu_N$-equivariant deformation of $Z$ with
quotient $\cX/S$. We have an isomorphism 
$$\omega_{\cZ/S} \otimes k(0) \cong \omega_Z \cong \cO_Z$$ 
by the base change property for the relative dualising sheaf.
Hence $\omega_{\cZ/S} \cong \cO_{\cZ}$ by Nakayama's lemma applied to the $\cO_{\cZ}$-module $\omega_{\cZ/S}$.
Thus  $\omega_{\cZ/S}^{\otimes N}$ is invertible and has a $\mu_N$-invariant generator.
Now, let $i \colon \cX^0 \inj \cX$ denote the Gorenstein locus of $\cX/S$ and $\pi^0  \colon \cZ^0 \map \cX^0$ 
the restriction
of the covering $\pi \colon \cZ \map \cX$. Then $\pi^0$ is an \'{e}tale $\mu_N$ quotient, hence
$$\omega_{\cX^0/S}^{\otimes N} \cong (\pi^0_{\star} \omega_{\cZ^0/S}^{\otimes N})^{\mu_N} \cong 
(\pi^0_{\star}\cO_{\cZ^0})^{\mu_N} \cong \cO_{\cX^0}.$$
Applying $i_{\star}$ we obtain $\omega_{\cX/S}^{[N]} \cong \cO_{\cX}$, thus $\cX/S$ is weakly $\bQ$-Gorenstein.
To prove that $N$ is the index, suppose $\omega_{\cX/S}^{[M]}$ is invertible for some $M \in \bN$,
and consider the natural map 
$$\omega_{\cX/S}^{[M]} \otimes k(0) \map \omega_X^{[M]}.$$
The map is an isomorphism in codimension $1$, and both sheaves are $S_2$, hence it is an isomorphism.
So $\omega_X^{[M]}$ is invertible and $N$ divides $M$. 
\end{proof}

\begin{lem} \label{lem-DVR}
Let $\cX/(0 \in T)$ be a flat family of slc surfaces over the germ of a curve.
Suppose that the general fibre is canonical, i.e., has only Du Val singularities,
and that $K_{\cX}$ is $\bQ$-Cartier. Then $\cX/T$ is $\bQ$-Gorenstein.
\end{lem}
\begin{proof}
We work locally at a point $P \in \cX$.
Let $Z \map X$ and $\cZ \map \cX$ be the canonical covering of $P \in X$ and $P \in \cX$ respectively.
Note that the index of $X$ equals the index of $\cX$ (\cite{KSB}, Lemma~3.16, p. 316), hence these maps have the same degree. 
We need to show that $\cZ/T$ is a deformation of $Z$.
Since the fibre $\cZ_0$ agrees with $Z$ over $X-P$, it is enough to show that $\cZ_0$ is Cohen-Macaulay.
The fibre $X$ of $\cX/T$ is slc, so the pair $(\cX,X)$ is log canonical --- this is an `inversion of adjunction' type result.
In more detail, after a finite surjective base change $T' \map T$, there is a semistable resolution
$\pi \colon \tilde{\cX'} \map \cX'= \cX \times_T T'$. 
Then the proof of \cite{KSB}, Theorem~5.1(a) shows that $\cX'/T'$ 
coincides with the canonical model of $\tilde{\cX'}$ over $\cX'$.
Hence $(\cX',X')/T'$ is log canonical. 
Finally, writing $g \colon \cX' \map \cX$ for the map induced by the base change $T' \map T$, we have 
$K_{\cX'} + X' = g^{\star}(K_{\cX}+X)$ by Riemann-Hurwitz,
so $(\cX,X)$ is log canonical by \cite{KM}, Proposition~5.20(4).
Since $X$ is Cartier and the general fibre is canonical it follows that $\cX$ is canonical.
Hence the cover $\cZ$ is also canonical, so in particular Cohen-Macaulay.
Then the  fibre $\cZ_0=(t=0) \subset \cZ$ is also Cohen-Macaulay.
\end{proof}

\subsection{Computing $\bQ$-Gorenstein deformations}

For $\cX/S$ a $\bQ$-Gorenstein family of slc surfaces, we define the canonical covering stack $\stX/S$ of the family
$\cX/S$, and show that the infinitesimal $\bQ$-Gorenstein deformations of $\cX/S$ correspond exactly to the infinitesimal 
deformations of $\stX/S$ (defined carefully below). 
We can then apply the results of \cite{I1},\cite{I2} to compute the $\bQ$-Gorenstein deformations of $\cX/S$ 
(Theorem~\ref{QGdefns}).
Note that, for our explicit computations in Sections~\ref{normal} and \ref{typeB}, we need only consider infinitesimal 
$\bQ$-Gorenstein deformations of an slc surface $X/\bC$. However, we must develop the theory 
for $\bQ$-Gorenstein families over an arbitrary affine scheme in order to 
establish `openness of versality' for $\bQ$-Gorenstein deformations (cf. \cite{Ar}, Section~4). 
This is used in the construction of the moduli space of stable pairs in Section~\ref{moduli}.

The following lemma motivates the definition of the canonical covering stack of a $\bQ$-Gorenstein family.

\begin{lem} \label{QGdefn_to_coverdefn}
Let $P \in X$ be an slc surface germ of index $N$ and $Z \map X$ the canonical covering with group $G \cong \mu_N$.
Let $\cZ/ (0 \in S)$ be a $G$-equivariant deformation of $Z$ inducing a $\bQ$-Gorenstein deformation 
$\cX/ (0 \in S)$ of $X$.
Then there is an isomorphism
$$\cZ \cong \underline{\Spec}_{\cX}(\cO_{\cX} \oplus \omega_{\cX/S} \oplus \cdots \oplus \omega^{[N-1]}_{\cX/S})$$
where the multiplication is given by fixing a trivialisation of $\omega_{\cX/S}^{[N]}$.
In particular, $\cZ/S$ is determined by $\cX/S$.
\end{lem}

\begin{proof}
Let $i \colon \cX^0 \inj \cX$ denote the open locus where the covering $\pi \colon \cZ \map \cX$ is \'{e}tale 
and let $\pi^0 \colon \cZ^0 \map \cX^0$ denote the restriction of the covering.
The map $\pi^0$ is an \'{e}tale $\mu_N$ quotient, hence
$$\cZ^0 \cong \underline{\Spec}_{\cX^0}(\cO_{\cX^0} \oplus \sL \oplus \cdots \oplus \sL^{\otimes N-1})$$
for some line bundle $\sL$ on $\cX^0$, with multiplication given by an isomorphism $\sL^{\otimes N} \cong \cO_{\cX}$.
Here the sheaves $\sL^{\otimes r}$ are the eigensheaves of the $G$ action on $\pi^0_{\star}\cO_{\cZ^0}$.
Since $\cZ$ is a deformation of the canonical covering of $Z$ of $X$, we may assume 
that the restriction of $\sL$ to the fibre $X^0$ is identified with $\omega_{X^0}$.
Now $\omega_{\cX^0/S}= (\pi^0_{\star}\omega_{\cZ^0/S})^G$ and $\omega_{\cZ/S} \cong \cO_{\cZ}$, hence
$\omega_{\cX^0/S}$ is isomorphic to a $G$-eigensheaf of $\pi^0_{\star}\cO_{\cZ^0}$ and so 
$\omega_{\cX^0/S} \cong \sL$ by our choice of $\sL$. 
Finally, $\cZ$ is determined by its restriction $\cZ^0$ since $\cZ$ is $S_2$ over $S$,
so we obtain an isomorphism as claimed.
\end{proof}

Let $\cX/S$ be a $\bQ$-Gorenstein family of slc surfaces. For $P \in \cX/S$ a point of index $N$, we define 
the \emph{canonical covering} $\pi \colon \cZ \map \cX$ of $P \in \cX/S$ by
$$\cZ=\underline{\Spec}_{\cX}(\cO_{\cX} \oplus \omega_{\cX/S} \oplus \cdots \oplus \omega^{[N-1]}_{\cX/S}),$$
where the multiplication is given by fixing a trivialisation of $\omega_{\cX/S}^{[N]}$ at $P$.
The canonical covering of $P \in \cX/S$ is uniquely determined in the \'{e}tale topology. 
Hence the data of canonical coverings 
everywhere locally on $\cX/S$ defines a Deligne-Mumford stack $\stX/S$ with coarse moduli space $\cX/S$,
the \emph{canonical covering stack} of $\cX/S$. 

The stack $\stX/S$ is flat over $S$ by Lemma~\ref{QGdefn_to_coverdefn}. 
Moreover, for any base change $T \map S$, let $\stX_T$ denote the canonical covering stack of
$\cX \times_S T/T$, then there is a canonical isomorphism $\stX_T \stackrel{\sim}{\map} \stX \times_S T$.
For, given an \'{e}tale neighbourhood $\cZ \map \stX$ as above, there is a corresponding \'{e}tale neighbourhood
$\cZ_T \map \stX_T$ and a natural map $\cZ_T \map \cZ \times_S T$ by the base change property for $\omega_{\cX/S}$.
The map is an isomorphism over the Gorenstein locus of $\cX \times_S T /T$ and both $\cZ_T$ and $\cZ \times_S T$ are
$S_2$ over $T$ by Lemma~\ref{QGdefn_to_coverdefn}, hence it is an isomorphism.

We collect some easy properties of the canonical covering stack $\stX/S$. 
There is a notion of an \'{e}tale map $U \map \stX$ and hence the notion of sheaves on the \'{e}tale site 
$\stX_{et}$ of the stack $\stX$. We shall only consider sheaves on $\stX_{et}$, and 
refer simply to `sheaves on $\stX$'. 
Let $\pi \colon \cZ \map \cX$ be a local canonical covering at $P \in \cX/S$, with group $G \cong \mu_N$. 
Then $\stX$ has local patch $[\cZ/G]$ over $P \in \cX$. 
Sheaves on $[\cZ/G]$ correspond to $G$-equivariant sheaves on $\cZ$.
Let $p \colon \stX \map \cX$ be the induced map to the coarse moduli space.
Thus, locally, $p$ is the map $[\cZ/G] \map \cZ/G$.
If $\cF$ is a sheaf on $[\cZ/G]$ and $\cF_{\cZ}$ is the corresponding $G$-equivariant sheaf on $\cZ$,
then $p_{\star}\cF=(\pi_{\star}\cF_{\cZ})^G$.
In particular, the functor $p_{\star}$ is exact. For, the map $\pi$ is finite and $(\pi_{\star}\cF_{\cZ})^G$
is a direct summand of $\pi_{\star}\cF_{\cZ}$ since we are in characteristic zero.

Let $A$ be a $\bC$-algebra and $A' \map A$ an infinitesimal extension.
Let $\cX/A$ be a $\bQ$-Gorenstein family of slc surfaces and $\stX/A$ the canonical covering stack of $\cX/A$.
A \emph{deformation} of $\stX/A$ over $A'$ is a Deligne-Mumford stack $\stX'/A'$, flat over $A'$, 
together with an isomorphism $\stX' \times_{\Spec A'} \Spec A \cong \stX$.
Observe that, since the extension $A' \map A$ is infinitesimal, we may identify the \'{e}tale sites
of $\stX'$ and $\stX$. Thus, equivalently, a deformation $\stX'/A'$ of $\stX/A$ is a sheaf $\cO_{\stX'}$ of flat $A'$-algebras
on the \'{e}tale site of $\stX$, together with an isomorphism $\cO_{\stX'} \otimes_{A'} A \cong \cO_{\stX}$.
From this point of view, infinitesimal deformations of stacks fit into the general framework of \cite{I1},\cite{I2}.
The stack $\stX/A$ is identified with the `ringed topos' over $A$ given by the \'{e}tale site of $\stX$ together 
with the structure sheaf $\cO_{\stX}$. The \emph{cotangent complex} $L_{\stX/A}$ of $\stX/A$ is a complex of $\cO_{\stX}$-modules
$L^i$ in degrees $i \le 0$, with $H^0(L_{\stX/A})=\Omega_{\stX/A}$. 
For an extension $A' \map A$  whose kernel $M$ satisfies $M^2=0$,
the groups $\Ext^i(L_{\stX/A},\cO_{\stX} \otimes_A M)$, $i=0,1,2$, control the deformations of $\stX/A$ over $A'$.
We refer to \cite{I2}, Section~1 for a review of cotangent complex theory, 
and to \cite{I1} for the definitive treatment.

In our calculations, we shall require the local-to-global spectral sequence for Ext and the Leray spectral sequence for
stacks. These are derived for ringed topoi, and thus for stacks, in \cite{SGA4}, Expos\'{e}~V.
In particular, if $\stX/A$ is the canonical covering stack of a $\bQ$-Gorenstein family  $\cX/A$
and $p \colon \stX \map \cX$ the induced map, then $H^i(\stX,\cF)=H^i(\cX,p_{\star}\cF)$ for $\cF$ a sheaf on $\stX$,
since $p_{\star}$ is exact.

\begin{notn}
Let $A$ be a $\bC$-algebra and $M$ a finite $A$-module.
For $\cX/A$ a flat family of schemes over $A$, let $L_{\cX/A}$ denote the cotangent complex of $\cX/A$. Define
$$T^i(\cX/A,M)=\Ext^i(L_{\cX/A}, \cO_{\cX} \otimes_A M)$$  
$$\cT^i(\cX/A,M)=\cExt^i(L_{\cX/A},\cO_{\cX} \otimes_A M)$$
For $\cX/A$ a $\bQ$-Gorenstein family of slc surfaces over $A$, let $\stX/A$ denote the canonical covering stack
of $\cX/A$ and $p \colon \stX \map \cX$ the induced map. Define
$$T^i_{QG}(\cX/A,M)=\Ext^i(L_{\stX/A}, \cO_{\cX} \otimes_A M)$$ 
$$\cT^i_{QG}(\cX/A,M)=p_{\star}\cExt^i(L_{\stX/A},\cO_{\stX} \otimes_A M)$$
\end{notn}

\begin{prop} \label{QG-stack}
Let $\cX/A$ be a $\bQ$-Gorenstein family of slc surfaces and $\stX/A$ the canonical covering stack.
Let $A' \map A$ be an infinitesimal extension of $A'$. For $\cX'/A'$ a $\bQ$-Gorenstein deformation
of $\cX/A$, let $\stX'/A'$ denote the canonical covering stack of $\cX'/A'$.
Then the map $\cX'/A' \mapsto \stX'/A'$ gives a bijection between the set of isomorphism classes of 
$\bQ$-Gorenstein deformations of $\cX/A$ over $A'$ and the set of isomorphism classes of deformations of $\stX/A$ over $A'$.
\end{prop}
\begin{proof}
If $\cX'/A'$ is a $\bQ$-Gorenstein deformation of $\cX/A$ then the canonical covering stack $\stX'/A'$
is a deformation of $\stX/A$.
Conversely, if $\stX'/A'$ is a deformation of $\stX/A$ then the coarse moduli space $\cX'/A'$ is a $\bQ$-Gorenstein
deformation of $\cX/A$.
It only remains to prove that, if $\stX'/A'$ is a deformation of $\stX/A$ with coarse moduli space $\cX'/A'$,
then the canonical covering stack $\tilde{\stX'}/A'$ of $\cX'/A'$ is isomorphic to $\stX'/A'$.
By induction, we may assume that the kernel $M$ of $A' \map A$ satisfies $M^2=0$.
Then the deformations of $\stX/A$ over $A'$ form an affine space under $T^1_{QG}(\cX/A,M)$ by \cite{I2}, Theorem~1.7.
Let $\stX'/A'$ and $\tilde{\stX'}/A'$ differ by an element $t \in T^1_{QG}(\cX/A,M)$; we show that $t=0$.
We have an exact sequence
$$0 \map H^1(\cT_{QG}^0(\cX/A,M)) \map T_{QG}^1(\cX/A,M) \stackrel{\theta}{\map} H^0(\cT_{QG}^1(\cX/A,M))$$ 
obtained from the local-to-global spectral sequence for $\Ext$ on the stack $\stX$.
The deformations $\tilde{\stX'}/A'$ and $\stX'/A'$ of $\stX/A$ induce isomorphic deformations
locally by Lemma~\ref{QGdefn_to_coverdefn}, hence $\theta(t)=0$, i.e., $t \in H^1(\cT_{QG}^0(\cX/A,M))$.
The natural map  $\cT_{QG}^0(\cX/A,M) \map \cT^0(\cX/A,M)$ is an isomorphism by Lemma~\ref{T^0}, 
so $t$ is identified with the element of $H^1(\cT^0(\cX/A,M))$ relating the deformations of $\cX/A$ induced by 
$\stX'/A'$ and $\tilde{\stX'}/A'$. But these deformations coincide by assumption, hence $t=0$ as required.
\end{proof}

\begin{lem} \label{T^0}
Let $\cX/A$ be a $\bQ$-Gorenstein family of slc surfaces and $M$ a finite $A$-module.
Then the natural map $\cT^0_{QG}(\cX/A,M) \map \cT^0(\cX/A,M)$ is an isomorphism.
\end{lem}
\begin{proof}
We work locally at $P \in \cX$.
Let $\pi \colon \cZ \map \cX$ be the canonical covering of $\cX/A$, with covering group $G$,
and $\stX=[\cZ/G]$ the canonical covering stack.
Then $\cT^0_{QG}(\cX/A,M)=(\pi_{\star}\cT^0(\cZ/A,M))^G$.
The natural map  $\cT^0_{QG}(\cX/A,M) \map \cT^0(\cX/A,M)$ is an isomorphism over the locus
where the covering $\pi$ is \'{e}tale, hence it suffices to show that   
$\cT^0_{QG}(\cX/A,M)$ and $\cT^0(\cX/A,M)$ are weakly $S_2$ over $A$.
First, we have
$$\cT^0(\cX/A,M)=\cHom(L_{\cX/A},\cO_{\cX} \otimes_A M)= \cHom(\Omega_{\cX/A},\cO_{\cX} \otimes_A M)$$
since the complex $L_{\cX/A}$ has cohomology $\Omega_{\cX/A}$ in degree $0$.
We claim that $\cO_{\cX} \otimes_A M$ is weakly $S_2$ over $A$, then $\cT^0(\cX/A,M)$ is weakly $S_2$ over $A$ by 
Lemma~\ref{S_2-basics}(1). To prove the claim, we may assume that $M=A/p$ for some prime ideal $p \subset A$ by
\ref{S_2-basics}(2). In this case $\cO_{\cX} \otimes_A M$ is $S_2$ over $A/p$ and so weakly $S_2$ over $A$ as desired. 
Second, the sheaf $\cT^0(\cZ/A,M)$ is weakly $S_2$ over $A$ as above, so $\pi_{\star}\cT^0(\cZ/A,M)$ is weakly $S_2$ over $A$.
Since $(\pi_{\star}\cT^0(\cZ/A,M))^G$ is a direct summand of $\pi_{\star}\cT^0(\cZ/A,M)$, it is also weakly $S_2$ over $A$.
\end{proof}

\begin{thm} \label{QGdefns}
Let $\cX_0/A_0$ be a $\bQ$-Gorenstein family of slc surfaces.
Let $M$ be a finite $A_0$-module.
\begin{enumerate}
\item The set of isomorphism classes of $\bQ$-Gorenstein deformations of $\cX_0/A_0$ over 
$A_0+M$ is naturally an $A_0$-module and is canonically isomorphic to  $T_{QG}^1(\cX_0/A_0,M)$.
Here $A_0+M$ denotes the ring $A_0[M]$, with $M^2=0$.
\item Let $A \map A_0$ be an infinitesimal extension and $A' \map A$ a further extension with kernel the 
$A_0$-module $M$. Let $\cX/A$ be a $\bQ$-Gorenstein deformation of $\cX_0/A_0$. 
\begin{enumerate}
\item There is a canonical element $o(\cX/A,A') \in T_{QG}^2(\cX_0/A_0,M)$ which vanishes if and only if
there exists a $\bQ$-Gorenstein deformation $\cX'/A'$ of $\cX/A$ over $A'$.
\item If $o(\cX/A,A')=0$, the set of isomorphism classes of $\bQ$-Gorenstein deformations $\cX'/A'$ is an affine space under
$T_{QG}^1(\cX_0/A_0,M)$.
\end{enumerate}
\end{enumerate}
\end{thm}
\begin{proof}
The $\bQ$-Gorenstein deformations of $\cX_0/A_0$ are identified with the deformations of the canonical covering stack 
$\stX_0/A_0$ of $\cX_0/A_0$ by Proposition~\ref{QG-stack}. Hence the theorem follows from \cite{I2}, 
Theorem~1.5.1 and Theorem~1.7. Note that, in part (2), we have used the
natural isomorphisms $T^i_{QG}(\cX_0/A_0,M) \stackrel{\sim}{\map} T^i_{QG}(\cX/A,M)$ given by \cite{I2}, 1.3.
\end{proof}

As remarked earlier, we need only consider infinitesimal deformations of an slc surface $X/\bC$ for our
later explicit computations. In the notation of the theorem, we may assume that $A_0=\bC$ and $M \cong \bC$.
We collect some useful notation and facts in this case below.
Define $T^i_X,\cT^i_X,T^i_{QG,X},\cT^i_{QG,X}$ by $T^i_{X}=T^i(X/\bC,\bC)$ etc. 
By the Theorem, first order $\bQ$-Gorenstein deformations of $X/\bC$ are identified with $T^1_{QG,X}$
and the obstructions to extending $\bQ$-Gorenstein deformations lie in $T^2_{QG,X}$.
We have $\cT^0_{QG,X}=\cT^0_X=\cHom(\Omega_X,\cO_X)$, the tangent sheaf of $X$, by Lemma~\ref{T^0}.
Working locally at $P \in X$, let $\pi \colon Z \map X$ be the canonical covering, with group $G$, then 
$\cT^i_{QG,X}=(\pi_{\star}\cT^i_Z)^G$. The sheaf $\cT^1_Z$ is supported on the singular locus of $Z$ and
$\cT^2_Z$ is supported on the locus where $Z$ is not a local complete intersection.
Finally, there is a local-to-global spectral sequence 
$$E^{pq}_2=H^p(\cT^q_X) \Rightarrow T^{p+q}_{QG,X}$$ 
given by the local-to-global spectral sequence for $\Ext$ on the canonical covering stack of $X$.

\subsection{Deformations of pairs}

Finally, we study deformations of stable pairs $(X,D)$.
We prove that the presence of the divisor $D$ does not produce any further obstructions.

\begin{defn}
Let $(P \in X,D)$ be a germ of a stable pair.
Let $N$ be the index of $X$ and $Z \map X$ the canonical covering, a $\mu_N$ quotient.
Let $D_Z$ denote the inverse image of $D$.
We say a deformation $(\cX,\sD)/(0 \in S)$ of $(X,D)$ is $\bQ$-\emph{Gorenstein}
if there is a $\mu_N$ equivariant deformation $(\cZ,\sD_{\cZ})/S$ of $(Z,D_Z)$
whose quotient is $(\cX,\sD)/S$.
\end{defn}
  
If $(\cX,\sD)/S$ is a $\bQ$-Gorenstein family of stable pairs
and $\pi \colon \cZ \map \cX/S$ is a local canonical covering of $\cX/S$, then
the closed subscheme $\sD_{\cZ} \inj \cZ$ is uniquely determined
by $\sD \inj \cX$. For, the ideal sheaf of $\sD_{\cZ}$
in $\cZ$ is $S_2$ over $S$ and agrees with the pullback of the 
ideal sheaf of $\sD$ in $\cX$ over the locus where $\pi$ is \'{e}tale.
Thus $\sD \inj \cX$ defines a closed substack $\stD \inj \stX$,
where $\stX$ is the canonical covering stack of $\cX/S$.

We first show that the families constructed in Theorem~\ref{thm-stable} satisfy the $\bQ$-Gorenstein condition.
This is needed to prove that the moduli space of stable pairs is proper.

\begin{lem}
Let $(\cX,\sD)/(0 \in T)$ be a flat family of stable pairs over the germ of a curve.
Suppose that the general fibre of $\cX/T$ is smooth
and that $K_{\cX}$ and $\sD$ are $\bQ$-Cartier. Then $(\cX,\sD)/T$ is $\bQ$-Gorenstein.
\end{lem}
\begin{proof}
The family $\cX/T$ is $\bQ$-Gorenstein by Lemma~\ref{lem-DVR}. Working locally at $P \in X \subset \cX$,
write 
\begin{eqnarray*}
\begin{array}{ccc}
(Z,D_Z) & \subset & (\cZ,\sD_{\cZ}) \\
\da &       & \da \\
(X,D) & \subset & (\cX,\sD) \\
\da &       & \da \\
0   & \in   & T
\end{array}
\end{eqnarray*}
for the canonical coverings together with the inverse images of the divisors $D$ and $\sD$.
We need to show that $\sD_{\cZ}$ is a deformation of $D_Z$.
We know that $\sD_{\cZ}$ is $\bQ$-Cartier and $D_Z$ is Cartier by Lemma~\ref{lem-DZCartier}; it follows that $\sD_{\cZ}$ is Cartier
(cf. \cite{KSB}, Lemma~3.16, p. 316) and thus $\sD_{\cZ} \otimes k(0)=D_Z$ as required.
\end{proof}

\begin{thm} \label{deformD}
Let $(\cX,\sD)/A$ be a $\bQ$-Gorenstein family of stable pairs.
Let $A' \map A$ be an infinitesimal extension and $\cX'/A'$ a $\bQ$-Gorenstein deformation of $\cX/A$
Then there exists a $\bQ$-Gorenstein deformation $(\cX',\sD')/A'$ of $(\cX,\sD)/A$.
\end{thm}
\begin{proof}
Let $\stX/A$ and $\stX'/A'$ denote the canonical covering stacks of $\cX/A$ and $\cX'/A'$,
and let $\stD \inj \stX$ be the closed substack determined by $\sD \inj \cX$.
We show that $\stD \inj \stX$ deforms to a closed substack $\stD' \inj \stX'$;
we then obtain the desired deformation $\sD' \inj \cX'$ of $\sD \inj \cX$
by forming the coarse moduli space.
By induction, we may assume that the kernel $M$ of $A' \map A$ satisfies $M^2=0$.
Then the obstruction to deforming  $\stD \inj \stX$ to 
a closed substack $\stD' \inj \stX'$ 
lies in $\Ext^2(L_{\stD / \stX},\cO_{\stD} \otimes_A M)$ by \cite{I2}, Theorem~1.7.
To complete the proof, we compute that this obstruction group is trivial.
The ideal sheaf $\cI$ of $\stD$ in $\stX$ is locally trivial,
i.e., $\stD$ is a Cartier divisor on $\stX$.
For, let $\cZ \map \cX$ be a  local canonical covering of $\cX/A$
and let $\sD_{\cZ} \inj \cZ$ be the closed subscheme corresponding
to $\stD \inj \stX$. Then $\sD_{\cZ}$ is flat over $A$ and has Cartier fibres 
by Lemma~\ref{lem-DZCartier}, hence $\sD_{\cZ}$ is Cartier.
In particular, the embedding $\stD \inj \stX$ is a local complete intersection,
thus $L_{\stD/\stX}$ is isomorphic to $\cI/\cI^2[-1]$ in the derived category of $\stD$, 
by \cite{I2}, p.160. Thus
$$\Ext^2(L_{\stD / \stX},\cO_{\stD} \otimes_A M) \cong \Ext^1(\cI/\cI^2,\cO_{\stD} \otimes_A M).$$
Now $\cI/\cI^2$ is locally trivial, hence $\cExt^1(\cI/\cI^2,\cO_{\stD} \otimes_A M)=0$ and 
$$\Ext^1(\cI/\cI^2,\cO_{\stD} \otimes_A M)=H^1(\stD,\cHom(\cI/\cI^2, \cO_{\stD} \otimes_A M)).$$
Next, we have
\begin{eqnarray*}
H^1(\stD,\cHom(\cI/\cI^2, \cO_{\stD} \otimes_A M)) & = & H^1(\stX,\cHom(\cI,\cO_{\stD}\otimes_A M)) \\
                                                   & = & H^1(\cX,p_{\star}\cHom(\cI,\cO_{\stD} \otimes_A M))
\end{eqnarray*}
where $p$ is the induced map $\stX \map \cX$. 
By cohomology and base change for $\cX/A$,  we may reduce to the case $A=M=\bC$; write $(X,D)=(\cX,\sD)$.
Applying $p_{\star}\cHom(\cI,-)$ to the exact sequence
$$0 \map \cI \map \cO_{\stX} \map \cO_{\stD} \map 0$$
of sheaves on $\stX$, we obtain the exact sequence
$$0 \map \cO_X \map \cO_X(D) \map p_{\star}\cHom(\cI,\cO_{\stD}) \map 0$$
of sheaves on $X$. Note that $\cHom(\cI,-)$ is exact since $\cI$ locally free,
and $p_{\star}$ is also exact. Consider the associated long exact sequence of cohomology
$$\cdots \map H^1(\cO_X(D)) \map H^1(p_{\star}\cHom(\cI,\cO_{\stD})) \map H^2(\cO_X) \map \cdots.$$
We have $H^1(\cO_X(D))=0$ by Lemma~\ref{lem-H1} and $H^2(\cO_X)=H^0(K_X)^{\vee}=0$ by Serre duality and
ampleness of $-K_X$. So $H^1(p_{\star}\cHom(\cI,\cO_{\stD}))=0$ as required.
\end{proof}

\begin{lem}\label{lem-DZCartier}
Let $(X,D)$ be a stable pair and $(Z,D_Z)$ a local canonical covering 
together with the inverse image of $D$. Then the divisor $D_Z$ is Cartier.
\end{lem}
\begin{proof}
If $d$ is divisible by $3$ then $\frac{d}{3}K_X+D \sim 0$, so $D_Z \sim -\frac{d}{3}K_Z \sim 0$.
Otherwise, by Theorem~\ref{thm-degnotmult3} and Propositions~\ref{prop-smoothableslt} and \ref{prop-singsP2}, 
the only possible singularities of $X$ are of the forms:
\begin{enumerate}
\item $\frac{1}{n^2}(1,na-1)$, where $3 \nd n$ and $(a,n)=1$.
\item $(xy=0) \subset \frac{1}{r}(1,-1,a)$, where $(a,r)=1$.
\end{enumerate}
In case (1), the local class group of $X$ is  $\bZ/n^2\bZ$.
So, since $dK_X+3D \sim  0$ and $3 \nd n$ , the divisor $D$ is locally a multiple of $K_X$, hence $D_Z$ is Cartier.
In case (2), the local class group of $\bQ$-Cartier divisors is $\bZ/r\bZ$,
generated by  $K_X$, so $D_Z$ is Cartier
\end{proof}

\begin{lem}\label{lem-H1}
Let $(X,D)$ be a stable pair. Then $H^1(\cO_X(D))=0$.
\end{lem}
\begin{proof}
We have $H^1(\cO_X(D))=H^1(\cO_X(K_X-D))^{\vee}$ by Serre duality and $-(K_X-D)$ is ample.
So if $X$ is log terminal our result follows by Kodaira vanishing.

Otherwise, let $\nu \colon X^{\nu} \map X$ be the normalisation of $X$ and $\tilde{\Delta} \map \Delta$ the normalisation
of the double curve $\Delta$. Then for any $\bQ$-Cartier divisor $E$ on $X$ there is an exact sequence
$$0 \map \cO_X(E) \map \cO_{X^{\nu}}(\nu^{\star}E) \map \cO_{\tilde{\Delta}}(\lfloor E|_{\tilde{\Delta}} \rfloor)$$
and hence a short exact sequence
$$0 \map \cO_X(E) \map \cO_{X^{\nu}}(\nu^{\star}E) \map \cF \map 0$$
where $\cF \inj \cO_{\tilde{\Delta}}(\lfloor E|_{\tilde{\Delta}} \rfloor)$.
Putting $E=K_X-D$, we have $H^0(\cF)=0$ since $-E$ is ample and $H^1(\cO_{X^{\nu}}(\nu^{\star}E))=0$ if 
$X^{\nu}$ is log terminal by Kodaira vanishing. So, in this case, the long exact sequence of cohomology gives
$H^1(\cO_X(K_X-D))=0$ as required. 

If $X^{\nu}$ is not log terminal, then $X$ is an elliptic cone by Theorems~\ref{thm-coarseclassn} and \ref{thm-ell_cone},
the degree $d$ is divisible by $3$ and 
$D \sim -\frac{d}{3}K_X$. An easy calculation shows that $H^1(\cO_X(D))=0$ in this case.
\end{proof}

\section{The moduli space of stable pairs} \label{moduli}

We construct the moduli space $\cM_d$ of stable pairs of degree $d$ using the deformation theory of Section~\ref{QG}.

\begin{defn} \label{defn-smoothable_deformation}
Let $(X,D)/\bC$ be a stable pair of degree $d$.
Let $(\cX^u,\sD^u) \map (0 \in S_0)$ be a versal $\bQ$-Gorenstein deformation of the pair $(X,D)/\bC$, 
where $S_0$ is of finite type over $\bC$. 
Let $S_1 \subset S_0$ be the open subscheme where the fibres of $\cX^u / S_0$ are isomorphic to $\bP^2$
and let $S_2$ be the scheme theoretic closure of $S_1$ in $S_0$.
A $\bQ$-Gorenstein deformation of $(X,D)$ is \emph{smoothable} if it is obtained by pullback from the deformation
$(\cX^u, \sD^u) \times_{S_0} S_2 \map (0 \in S_2)$.
\end{defn}

\begin{rem}
This definition is vacuous if the degree is not a multiple of $3$, i.e., any $\bQ$-Gorenstein deformation of $(X,D)$ 
is automatically smoothable.
For $(X,D)$ has unobstructed $\bQ$-Gorenstein deformations if $3 \nd d$ by Theorem~\ref{thm-Msmooth}, so that, in the notation above, 
the germ $0 \in S_0$ is smooth and thus the open subscheme $S_1 \subset S_0$ is dense and $S_2=S_0$. 
However, if $d$ is a multiple of $3$, there are examples where $S_0$ is reducible and $S_2$ is an irreducible component of $S_0$.
\end{rem}

\begin{defn}
Let $\underline{Sch}$ be the category of noetherian schemes over $\bC$. Let $d \in \bN$, $d \ge 4$.  
We define a stack $\cM_d \map \underline{Sch}$ as follows:
$$\cM_d(S) = \left\{ (\cX,\sD)/S \left| \begin{array}{c}        
                                        \mbox{$(\cX,\sD)/S$ is a $\bQ$-Gorenstein smoothable}\\
                                        \mbox{family of stable pairs of degree $d$}
                                        \end{array} \right. \right\}$$
\end{defn}

\begin{thm}
The stack $\cM_d$ is a separated and proper Deligne--Mumford stack.
The underlying coarse moduli space is a compactification 
of the moduli space of smooth plane curves of degree $d$.
\end{thm}

We give the salient points in the proof of the theorem.
Using the obstruction theory for $\bQ$-Gorenstein deformations obtained in Section~\ref{QG},
we deduce the existence of versal $\bQ$-Gorenstein deformations for stable pairs, corresponding to local
patches of the stack $\cM_d$ \cite{Ar}. 
To prove boundedness, i.e., that only finitely many patches are required, we first 
bound the index of a surface $X$ occurring in a stable pair of degree $d$ (Theorem~\ref{thm-index}).
Then, letting $N(d) \in \bN$ be such that $N(d)K_X$ is Cartier for each such pair $(X,D)$,
we have a polarisation on each $(X,D)$ given by $-N(d)K_X$. The Hilbert polynomial
is fixed by our smoothability assumption; hence boundedness follows by \cite{Ko}, Theorem~2.1.2.
Finally, the stack $\cM_d$ is separated and proper by Theorem~\ref{thm-stable}.

\begin{thm} \label{thm-index}
Let $(X,D)$ be a stable pair of degree $d$.
Then the index of each point $P \in X$ is at most $d$.
Moreover, the same result holds if $(X,D)$ is a semistable pair of degree $d$
and $d$ is not a multiple of $3$. 
\end{thm}
\begin{proof}
The pair $(X,(\frac{3}{d}+\epsilon)D)$ is slc, hence $D$ misses the strictly slc points of $X$.
Then the condition $dK_X+3D \sim 0$ shows that the index of $X$ is at most $d$ at such points.

The slt singularities of $X$ are of the following types (by Propositions~\ref{prop-smoothableslt} and \ref{prop-singsP2}):
\begin{enumerate}
\item $\frac{1}{n^2}(1,na-1)$, where $(a,n)=1$ and $3 \nd n$.
\item $(xy=0) \subset \frac{1}{r}(1,-1,a)$ where $(a,r)=1$.
\item $(x^2=zy^2) \subset \bA^3$.
\end{enumerate}
The index of $X$ equals $n$, $r$ and $1$ in cases $(1)$, $(2)$ and $(3)$ respectively.

In case (1) let $\tilde{X} \map X$ be the local smooth covering of $X$ and $\tilde{D}$ the inverse image of $D$.
Write $\tilde{X}=\bA^2_{x,y}$ and $\tilde{D}=(f(x,y)=0)$.
The multiplicity of the divisor $\tilde{D}$ at $0 \in \tilde{X}$ is 
strictly less than $\frac{2d}{3}$ since $(\tilde{X},(\frac{3}{d}+\epsilon)\tilde{D})$ is log canonical. 
Let $x^iy^j$ be a monomial appearing in the polynomial $f(x,y)$ such that $i+j$ is minimal,
thus $i+j < \frac{2d}{3}$. Then $3(i+(na-1)j)=dna \mod n^2$, using $dK_X+3D \sim 0$.
In particular, $i=j \mod n$. Thus if $n > d$ then $i=j< \frac{d}{3}$ and $3i=d \mod n$, a contradiction. 

In case (2) let $\tilde{X} \map X$ be the canonical covering of $X$, let $\tilde{D}$ denote the inverse image of $D$
and $\tilde{\Delta}$ the inverse image of the double curve of $X$.
Write $\tilde{X}=(xy=0) \subset \bA^3_{x,y,z}$ and $\tilde{D}=(f(x,y,z)=0)$.
Then $$\tilde{D}|_{\tilde{\Delta}}=(f(0,0,z)=0)=(z^k+\cdots=0) \subset \bA^1_z,$$  
where $k$ is the multiplicity of $\tilde{D}|_{\Delta}$ at $0 \in \tilde{\Delta}$.
Then $k < \frac{d}{3}$ since $(X,(\frac{3}{d}+\epsilon)D)$ is slc and $3k=d \mod r$ since $dK_X+3D \sim 0$.
Hence $r \le d$ as required.

If $(X,D)$ is a semistable pair of degree $d$, then $X$  has only singularities of type $\frac{1}{n^2}(1,na-1)$, the pair 
$(X,\frac{3}{d}D)$ is log canonical and $dK_X+3D \sim 0$.
Then, assuming $3 \nd d$, proceeding as in case (1) above we deduce the same result.
\end{proof}

\section{A coarse classification of the degenerate surfaces} \label{coarse}

If $(X,D)$ is a stable pair then the surface $X$ is slc and the divisor $-K_X$ is ample. 
We use these two properties to obtain a coarse classification of the possible surfaces $X$. 
We first describe the pairs $(Y,C)$ where $Y$ is an irreducible component of the normalisation of $X$ 
and $C$ is the inverse image of the double curve of $X$. We then glue such pairs together to obtain the 
classification of the surfaces $X$.

\begin{thm} (\cite{KM}, p. 119, Theorem~4.15) \label{thm-lcsings}
Let $P \in Y$ be the germ of a surface and $C$ an effective divisor on  $Y$ such that the pair $(Y,C)$ is log canonical.
Then, assuming $C \neq 0$, the germ $(P \in Y,C)$ is of one of the following types:
\begin{enumerate}
\item $(\frac{1}{r}(1,a),(x=0))$, where $(a,r)=1$.
\item $(\frac{1}{r}(1,a),(xy=0))$, where $(a,r)=1$.
\item $(\frac{1}{r}(1,a),(xy=0))/ \mu_2$, where the $\mu_2$ action is etale in codimension $1$ and interchanges 
$(x=0)$ and $(y=0)$.
\end{enumerate}
Moreover (1) is log terminal, whereas (2) and (3) are strictly log canonical.
\end{thm}

\begin{notn}
We denote singularities of types (1), (2) and (3) by 
$(\frac{1}{r}(1,a),\Delta)$, $(\frac{1}{r}(1,a),2\Delta)$ and $(D,\Delta)$ respectively. 
The $D$ stands for dihedral --- the surface singularities $P \in Y$ here include the dihedral Du Val singularities.
\end{notn}

\begin{thm}\label{thm-cpts}
Let $Y$ be a surface and $C$ an effective divisor on $Y$ such that the pair $(Y,C)$ is log canonical
and $-(K_Y+C)$ is ample. Then $(Y,C)$ is of one of the following types:
\begin{enumerate}
\item[(I)] $C=0$.
\item[(II)] $C \cong \bP^1$ and $(Y,C)$ is log terminal.
\item[(III)] $C \cong \bP^1 \cup \bP^1$, where the components meet in a single node.
\item[(IV)] $C \cong \bP^1$ and $(Y,C)$ has a singularity of type $(D,\Delta)$.
\end{enumerate}
Moreover, in case (I) the surface $Y$ has at most one strictly log canonical singularity, 
in case (III) the pair $(Y,C)$ is log terminal away from the node of $C$ 
and in case (IV) the pair $(Y,C)$ is log terminal away from the singularity of type $(D,\Delta)$.
\end{thm}

\begin{proof}
The pair $(Y,C)$ is log canonical and $-(K_Y+C)$ is ample by assumption, hence the locus where $(Y,C)$ is not klt is 
connected by the connectedness theorem of Koll\'{a}r and Shokurov (cf. \cite{KM}, p. 173, Theorem~5.48 and Corollary~5.49).
In other words, either $C=0$ and $Y$ has at most one strictly log canonical singularity, or $C$ is connected and $Y$ is log 
terminal away from $C$. 

If $C \neq 0$, let $\Gamma$ be a component of $C$. Then
$$(K_Y+C)\Gamma = (K_Y+ \Gamma)\Gamma + (C- \Gamma)\Gamma = 2p_a(\Gamma)-2 + \Diff(Y,\Gamma) + (C- \Gamma)\Gamma$$
where $\Diff(Y,\Gamma)$ is the different of the pair $(Y,\Gamma)$, i.e., the correction to the adjunction formula for 
$\Gamma \subset Y$ required due to the singularities of $Y$ at $\Gamma$ (\cite{FA}, Chapter~16).
Now $(K_Y+C)\Gamma <0$ since $-(K_Y+C)$ is ample, the different $\Diff(Y,\Gamma) \ge 0$ and $(C-\Gamma)\Gamma \ge 0$.
So $p_a(\Gamma)=0$, i.e., the curve $\Gamma$ is smooth and rational, and 
$$\Diff(Y,\Gamma) + (C- \Gamma)\Gamma <2.$$
The singularities of $(Y,C)$ at $\Gamma$ are of the forms
$(\frac{1}{r}(1,a),\Delta)$, $(\frac{1}{r}(1,a),2\Delta)$ and $(D,\Delta)$ as described in Theorem~\ref{thm-lcsings}.
We calculate that these singularities contribute $1-\frac{1}{r}$, $1$ and $1$ to the value of 
$\Diff(Y,\Gamma) + (C- \Gamma)\Gamma$ respectively.
The theorem now follows easily.
\end{proof}

\begin{notn}
Let $X$ be an slc surface.
Let $\Delta$ denote the double curve of $X$. 
Let $X_1,\ldots, X_n$ be the irreducible components of $X$ and write $\Delta_i$ for the restriction of $\Delta$ to $X_i$.
Let $\nu \colon X^{\nu} \map X$ be the normalisation of $X$ and $\Delta^{\nu}$ the inverse image of $\Delta$;
also write $X_i^{\nu}$ for the normalisation of $X_i$ and $\Delta^{\nu}_i$ for the inverse image of $\Delta_i$.

The map $\Delta^{\nu} \map \Delta$ is 2-to-1.
Let $\Gamma \subset \Delta$ be a component and write $\Gamma^{\nu}$ for its inverse image on $X^{\nu}$.
Then either $\Gamma^{\nu}$ has two components mapping birationally to $\Gamma$ or $\Gamma^{\nu}$ is irreducible and is a
double cover of $\Gamma$.
In the latter case we say that the curve $\Gamma^{\nu} \subset X^{\nu}$ is \emph{folded} to obtain $\Gamma \subset X$. 
\end{notn} 

\begin{thm}\label{thm-coarseclassn}
Let $X$ be a slc surface such that $-K_X$ is ample.
Then $X$ is of one of the following types:
\begin{enumerate}
\item[(A)] $X$ is normal. 
\item[(B)] $X$ has two normal components meeting in a smooth rational curve and is slt.
\item[(B*)] $X$ is irreducible, non-normal and slt. The pair $(X^{\nu},\Delta^{\nu})$ is of type II and 
$X$ is obtained by
folding the curve $\Delta^{\nu}$.
\item[(C)] $X$ has $n$ components $X_1, \cdots, X_n$ such that $(X^{\nu}_i,\Delta^{\nu}_i)$ is of type III for each $i$.
One component of $\Delta^{\nu}_i$ is glued to a component of $\Delta^{\nu}_{i+1 \mod n}$ for each $i$ so that the nodes 
of the curves $\Delta^{\nu}_i$ coincide and the components $X_i$ of $X$ form an `umbrella'.
\item[(D)] $X$ has $n$ components $X_1,\ldots, X_n$ such that $(X_i,\Delta_i)$ is of type III for $2 \le i \le n-1$.
Either $(X_1,\Delta_1)$ is of type IV or 
$(X_1^{\nu},\Delta^{\nu}_1)$ is of type III and $(X_1,\Delta_1)$ is obtained by folding one component of 
$\Delta^{\nu}_1$; similiarly for $(X_n,\Delta_n)$. The components $(X_1,\Delta_1), \ldots ,(X_n,\Delta_n)$ are glued
sequentially so that the nodes of the curves $\Delta_i$ and any $(D,\Delta)$ singularities on $(X_1,\Delta_1)$ and 
$(X_n,\Delta_n)$ coincide and the components $X_i$ of $X$ form a `fan'.  
\end{enumerate}
\end{thm}

\begin{proof}
Let $(Y,C)$ be a component of the pair $(X^{\nu},\Delta^{\nu})$, then $(Y,C)$ is log canonical and $-(K_Y+C)$ is ample
since $X$ is slc and $-K_X$ is ample. Hence $(Y,C)$ is of one of the types I--IV described in Theorem~\ref{thm-cpts}.
Glueing the components back together (using the classification of slc singularities \cite{KSB}) we obtain the 
classification of the surfaces $X$ given above.
\end{proof}

\begin{thm} \label{thm-rationality}
Let $(Y,C)$ be as in Theorem~\ref{thm-cpts}.
Then either $Y$ is rational or $C=0$ and $Y$ is an elliptic cone.
\end{thm}
\begin{proof}
Let $\pi \colon \tilde{Y} \map Y$ be the minimal resolution of $Y$.
Let $\tilde{C}$ be the $\bQ$-divisor defined by the equation
$$K_{\tilde{Y}}+\tilde{C}=\pi^{\star}(K_Y+C).$$
Note that $\tilde{C}$ is effective since $\pi$ is minimal and $-(K_{\tilde{Y}}+\tilde{C})$ is nef and big since 
$-(K_Y+C)$ is ample. So, running a MMP, 
we obtain a birational morphism $\phi \colon \tilde{Y} \map Y_1$ where either $Y_1 \cong \bP^2$ or $Y_1$ has the
structure of a $\bP^1$-bundle $q \colon Y_1 \map B$ over a smooth curve $B$.
We may assume that $Y$ is not rational, so that we are in the second case and the curve $B$ has positive genus.

We claim that there is an irrational component of the divisor $\tilde{C}$. For, otherwise, the image $C_1=\phi_{\star}\tilde{C}$ 
is a sum of 
fibres of the ruling. Then $-(K_{Y_1}+C_1)$ nef and big implies that $-K_{Y_1}$ is nef and big, hence
$h^1(\cO_{Y_1})=0$ by Kodaira vanishing. So $B$ has genus zero, a contradiction.

We have $\Supp \tilde{C} \subset C' \cup \Exc(\pi)$ where $C'$ denotes the strict transform of $C$ on $\tilde{Y}$
and $\Exc(\pi)$ is the exceptional locus of $\pi$. 
By Theorem~\ref{thm-cpts} the curve $C$ has only rational components, so $\Exc(\pi)$ contains an irrational curve and
$Y$ has a simple elliptic singularity by the classification of log canonical singularities. 
Let $E$ denote the corresponding
$\pi$-exceptional elliptic curve on $\tilde{Y}$. Then $E$ has multiplicity $1$ in $\tilde{C}$ and is
horizontal with respect to the birational ruling $\tilde{Y} \map B$. 
The divisor $-(K_{\tilde{Y}}+\tilde{C})$ is big, so $-(K_{\tilde{Y}}+\tilde{C})f >0$
where $f$ is a fibre of the ruling. Hence $E \cdot f=1$, i.e., the curve $E$ is a section of the ruling.

We show that $\tilde{Y}$ is actually biregularly ruled over the elliptic curve $E$.
Suppose not, then there is a degenerate fibre; write $A$ for a component meeting $E$.
Then $A$ is not contained in $\Supp \tilde{C}$ and $(K_{\tilde{Y}}+\tilde{C})A \le 0$, with equality if and only if
$A$ is contracted by $\pi$. But also
$$(K_{\tilde{Y}}+\tilde{C})A \ge K_{\tilde{Y}}A+E \cdot A \ge -1 +1=0$$
with equality only if $A$ is a $(-1)$-curve. So $A$ is a $(-1)$-curve which is contracted by $\pi$, a contradiction since
$\pi$ is minimal. 

Thus $\tilde{Y}$ is a $\bP^1$-bundle over an elliptic curve and the surface $Y$ is obtained by 
contracting the negative section; so $Y$ is an elliptic cone. Finally $C=0$ by Theorem~\ref{thm-cpts}.
\end{proof}

\section{Preliminary smoothability results} \label{prelim}

We collect some further restrictions on the degenerate surfaces $X$ 
implied by the existence of a smoothing to the plane.
We give a more detailed analysis of the possible singularities and some restrictions 
on the Picard numbers of the components of $X$. We deduce that a surface of type B* 
cannot smooth to the plane.

\begin{prop} (\cite{KSB}, Theorem~4.23 and 5.2) \label{prop-smoothableslt}
Let $P \in X$ be an slt surface singularity which admits a $\bQ$-Gorenstein smoothing.
Then $P \in X$ is of one of the following types:
\begin{enumerate}
\item A Du Val singularity.
\item $\frac{1}{dn^2}(1,dna-1)$, where $(a,n)=1$.
\item $(xy=0) \subset \frac{1}{r}(1,-1,a)$, where $(a,r)=1$.
\item $(x^2=y^2z) \subset \bA^3$.
\end{enumerate}
\end{prop}

\begin{prop} \label{prop-singsP2} 
Let $X$ be an slc surface which admits a $\bQ$-Gorenstein smoothing to $\bP^2$.
Then the log terminal singularities of $X$ are of the form \mbox{$\frac{1}{n^2}(1,na-1)$}, where $3 \nd n$.  
\end{prop}

\begin{proof}
In the case that $X$ is globally log terminal this was proved in \cite{Ma}, Section~3.
The same argument proves our result.
\end{proof}

\begin{prop} \label{prop-rho}
Let $X$ be an slc surface such that $-K_X$ is ample.
Suppose that $X$ admits a smoothing to $\bP^2$, i.e., there exists a flat family \mbox{$\cX/(0 \in T)$} over the germ of a curve with
special fibre $X$ and general fibre $\bP^2$. Then, in the cases A, B and B* of Theorem~\ref{thm-coarseclassn} the Picard numbers 
of the components of $X$ are as follows: 
\begin{enumerate}
\item[(A)] $\rho(X)=1$.
\item[(B)] Either (i) $\rho(X_1)=\rho(X_2)=1$ or (ii) $\{\rho(X_1),\rho(X_2)\}=\{1,2\}$.
\item[(B*)] $\rho(X^{\nu})=1$. 
\end{enumerate}
Moreover, given a smoothing $\cX/T$ of $X$ as above, 
the total space $\cX$ is $\bQ$-factorial unless $X$ is of type B, case (i).
\end{prop}

\begin{rem} 
In fact a surface of type B* never admits a smoothing to $\bP^2$  by Theorem~\ref{thm-B*} below
--- the above result is required in the proof.
\end{rem} 

\begin{proof}
Consider the commutative diagram
\begin{eqnarray*}
\begin{array}{ccc}
\Pic \cX & \map & \Pic X \\ 
\da    &      &   \da \\
H^2(\cX,\bZ) & \map & H^2(X,\bZ)
\end{array}
\end{eqnarray*} 
We claim that all these maps are isomorphisms.
First, the restriction map $H^2(\cX,\bZ) \map H^2(X,\bZ)$ is an isomorphism because $X$ is a homotopy retract of $\cX$.
Second, the map $c_1 \colon \Pic X \map H^2(X,\bZ)$ fits into the long exact sequence of cohomology
$$\cdots \map H^1(\cO_X) \map \Pic X \map H^2(X,\bZ) \map H^2(\cO_X) \map \cdots$$
associated to the exponential sequence on $X$. Now $H^2(\cO_X)=0$ since $-K_X$ is ample, so also $H^1(\cO_X)=0$  
since $\chi(\cO_X)=\chi(\cO_{\bP^2})$. Thus $c_1$ is an isomorphism as claimed. Similiarly
$c_1 \colon \Pic \cX \map H^2(\cX,\bZ)$ is also an isomorphism.
Hence the restriction map $\Pic(\cX) \map \Pic(X)$ is an isomorphism. 
Now $\Cl(\cX) \cong \bZ^n$ by Lemma~\ref{classgp}, so we have the inequality
$$\rho(X) = \dim \Pic(X) \otimes \bQ = \dim \Pic(\cX) \otimes \bQ \le \dim \Cl(\cX) \otimes \bQ = n,$$
with equality if and only if $\cX$ is $\bQ$-factorial.

We next relate the Picard numbers of $X$ and its irreducible components.
Note that $\rho(X) =\dim H_2(X,\bQ)$ by the above and similiarly for the components of $X$, since they are rational by 
Theorem~\ref{thm-rationality}.  
If $X$ is a surface of type~B then the Mayer-Vietoris sequence for $X=X_1 \cup X_2$ yields the following exact sequence:
\begin{eqnarray*}
\begin{array}{ccccccccc}
0 & \map & \bQ & \map & H_2(X_1,\bQ)  \oplus  H_2(X_2,\bQ) & \map & H_2(X) & \map & 0\\
  &      &   1 & \mapsto &   \Delta    \oplus     \Delta      &      &        &      &
\end{array}
\end{eqnarray*}
So $\rho(X)=\rho(X_1)+\rho(X_2)-1$.
If $X$ is a surface of type B* then an easy Mayer-Vietoris argument shows that $H_2(X^{\nu}) \cong H_2(X)$, so $\rho(X^{\nu})=\rho(X)$.
Our result now follows easily using the inequality $\rho(X) \le n$ derived above.
\end{proof}

\begin{thm} \label{thm-B*}
Let $X$ be a surface of type B*. Then $X$ does not admit a smoothing to $\bP^2$.
\end{thm}
\begin{proof}
Suppose $X$ is a counter-example and let $\cX/T$ be a smoothing of $X$ to $\bP^2$. 
Then $\rho(X^{\nu})=1$ and $\cX$ is $\bQ$-factorial by Proposition~\ref{prop-rho};
in particular $K_{\cX}$ is $\bQ$-Cartier.
Thus
$$(K_{X^{\nu}}+\Delta^{\nu})^2=(\nu^{\star}K_X)^2=K_X^2=K_{\bP^2}^2=9$$
and so $K_{X^{\nu}}^2 > 9$, since $-K_X$ is ample and $\rho(X^{\nu})=1$.
Hence $$K_{X^{\nu}}^2+\rho(X^{\nu}) > 10.$$ 
On the other hand, let $\tilde{X}$ be the minimal resolution of $X^{\nu}$, then 
\mbox{$K_{\tilde{X}}^2+\rho(\tilde{X})=10$} by Noether's formula, since $\tilde{X}$ is rational by Theorem~\ref{thm-rationality}.
We calculate below that the only possible singularities on $X^{\nu}$ cause an increase in $K^2+\rho$ 
on passing to the minimal resolution, so we have a contradiction. In fact, if a \emph{normal} rational surface singularity
$P \in S$ admits a $\bQ$-Gorenstein smoothing then, writing $\tilde{S} \map S$ for minimal resolution, 
the Milnor number $\mu$ of the smoothing equals $K_{\tilde{S}}^2+\rho(\tilde{S})$ (cf. \cite{Lo}).
In particular, $K_{\tilde{S}}^2+\rho(\tilde{S})$ is a non-negative integer.

The pair $(X^{\nu},\Delta^{\nu})$ has singularities of types $(\frac{1}{n^2}(1,na-1),0)$ and $(\frac{1}{r}(1,a),\Delta)$, 
with the latter cases occurring in pairs $\frac{1}{r}(1,a)$ and $\frac{1}{r}(-1,a)$, by Proposition~\ref{prop-smoothableslt}
and Proposition~\ref{prop-singsP2}. 
Given a cyclic quotient singularity $\frac{1}{r}(1,a)$, let 
$\frac{r}{a}=[b_1,\ldots,b_k]$ be the expansion of $\frac{r}{a}$ as a Hirzebruch--Jung continued fraction 
(\cite{Fu}, pp. 45-7). Then the minimal resolution of the singularity has exceptional locus
a chain of smooth rational curves $E_1,\ldots,E_n$ with self intersections $-b_1,\ldots,-b_k$. 
On passing to the minimal resolution, the change in $K^2+\rho$ is given by
$$\delta= E^2+4-\frac{1}{r}(a+a'+2)$$
where $E=E_1+\cdots +E_n$ and $a'$ is the inverse of $a$ modulo $r$.
For the singularity $\frac{1}{n^2}(1,na-1)$ we calculate $\delta=0$ using the inductive description of the minimal 
resolutions of these singularities (see \cite{KSB}, p. 314, Proposition~3.11).
For a pair of singularities $\frac{1}{r}(1,a)$, $\frac{1}{r}(1,-a)$ we calculate $\delta_1+\delta_2=4(1-\frac{1}{r})$.
Here we use the following elementary fact: if $\frac{r}{a}=[b_1,\ldots,b_k]$ and $\frac{r}{r-a}=[c_1,\ldots,c_l]$, then
$$\sum (b_i-1) = \sum (c_j-1) =k+l-1.$$
\end{proof}

\section{Simplifications in the case $3 \nd d$} \label{simp}

We show how the classification of stable pairs of degree $d$ simplifies considerably if $d$ is not a multiple of $3$.
We deduce that the stack $\cM_d$ is smooth in this case.

\begin{thm}\label{thm-degnotmult3}
Let $(X,D)$ be a stable pair of degree $d$, where $d$ is not a multiple of $3$.
Then $X$ is slt. So $X$ is either a normal log terminal surface or a surface of type $B$.
In particular, the surface $X$ has either $1$ or $2$ components.
\end{thm}

\begin{proof}
Since $3 \nd d$, the condition $dK_X+3D \sim 0$ gives an arithmetic condition on $X$, namely that $K_X$ is $3$-divisible as a 
$\bQ$-Cartier divisor class on $X$. Roughly, given a curve $\Gamma$ on $X$ we have $dK_X \cdot \Gamma= 3D \cdot \Gamma$, 
so $K_X \cdot \Gamma$ should be divisible by $3$. Of course the intersection number $D \cdot \Gamma$ is not an integer in general
since $X$ is singular, but this is the idea of the proof.

First suppose that $X$ is normal. Then either $X$ is log terminal or $X$ is an elliptic cone by Theorem~\ref{thm-ell_cone}.
In the second case, let $\Gamma$ be a ruling of the cone, then $3D \cdot \Gamma =-dK_X \cdot \Gamma=d$.
But $D$ misses the singularity of $X$ since the pair $(X,(\frac{3}{d}+\epsilon)D)$ is log canonical, hence $D$ is Cartier
and $D \cdot \Gamma$ is an integer, a contradiction.

Now suppose $X$ is not normal. 
Let $(Y,C)$ be a component of the pair $(X^{\nu},\Delta^{\nu})$, where $X^{\nu}$ is the normalisation of $X$ and $\Delta^{\nu}$ is
the inverse image of the double curve of $X$. We need to show that $(Y,C)$ is log terminal. 
Suppose this is not the case, then the pair $(Y,C)$ has a singularity of type $(\frac{1}{r}(1,a),2\Delta)$ or $(D,\Delta)$.
Let $\Gamma$ be a component of $C$ passing through such a point. Then there is at most one other singularity of 
$(Y,C)$ on $\Gamma$, of type $(\frac{1}{s}(1,b),\Delta)$, and
$$(K_Y+C)\Gamma=-2+1+(1-\frac{1}{s})=-\frac{1}{s},$$ 
cf. Theorem~\ref{thm-cpts}. We allow $s=1$, corresponding to the case that there are no further singularities at $\Gamma$.
Then $3D \cdot \Gamma = -d(K_Y+C)\Gamma=\frac{d}{s}$. But $D$ misses the strictly log canonical singularity of $X$,
hence $sD$ is Cartier near $\Gamma$ and $sD \cdot \Gamma$ is an integer, a contradiction.
\end{proof}

\begin{thm}\label{thm-Msmooth}
The stack $\cM_d$ is smooth if $3 \nd d$.
\end{thm}
\begin{proof}
Assume $3 \nd d$.
Given a stable pair $(X,D)$ of degree $d$, the surface $X$ is slt by Theorem~\ref{thm-degnotmult3}.
So $X$ is either normal and log terminal or a surface of type B, and $X$ has unobstructed $\bQ$-Gorenstein deformations by
Theorem~\ref{thm-manetti_defns} or Theorem~\ref{thm-B} respectively.
The $\bQ$-Gorenstein deformations of the pair $(X,D)$ are thus unobstructed by Theorem~\ref{deformD}. Hence 
$\cM_d$ is smooth as required.
\end{proof}

\section{The normal surfaces} \label{normal}
\subsection{Log terminal surfaces}
Log terminal degenerations of the plane have been classified by Manetti \cite{Ma}.
We announce a refinement of his result below (Theorem~\ref{thm-manetti}), the proof will appear elsewhere.

\begin{defn} A \emph{Manetti surface} is a normal log terminal surface which smoothes to $\bP^2$.
\end{defn}

\begin{thm} \label{thm-manetti_defns}
Let $X$ be a Manetti surface. Then $X$ has unobstructed $\bQ$-Gorenstein deformations.
\end{thm}
\begin{proof}
The obstructions are contained in $T^2_{QG,X}$ and there is a spectral sequence
$$E_2^{pq} =H^p(\cT^q_{QG,X}) \Rightarrow T^{p+q}_{QG,X},$$
hence it is enough to show that $H^p(\cT_{QG,X}^q)=0$ for $p+q=2$.
The sheaf $\cT^1_{QG,X}$ is supported on the singular locus, a finite set, so $H^1(\cT^1_{QG,X})=0$.
The singularities of $X$ are of the form \mbox{$\frac{1}{n^2}(1,na-1)$} by Proposition~\ref{prop-singsP2}.
Let $\pi \colon Z \map X$ be a local canonical covering, with group $\mu_n$. 
Then $Z$ is a hypersurface, so $\cT^2_Z=0$ and $\cT^2_{QG,X}=(\pi_{\star}\cT^2_Z)^{\mu_n}=0$.
Hence $H^0(\cT^2_{QG,X})=0$. Finally $H^2(\cT^0_{QG,X})=H^2(\cT_X)=0$ by \cite{Ma}, p. 113.
\end{proof}

\begin{thm} \label{thm-manetti}
Let $(a,b,c)$ be a solution of the Markov equation $$a^2+b^2+c^2=3abc.$$
Then the weighted projective space $X=\bP(a^2,b^2,c^2)$ smoothes to $\bP^2$.
Moreover, $X$ has no locally trivial deformations and there is precisely one 
deformation parameter for each singularity.
Conversely, every Manetti surface is obtained as a $\bQ$-Gorenstein deformation of such a surface $X$.
\end{thm}

The solutions of the Markov Equation are easily described \cite{Mo}. First, $(1,1,1)$ is a solution. Second, given one solution, we obtain
another by regarding the equation as a quadratic in one of the variables, $c$ (say), and replacing $c$ by the other root, i.e.,
$$(a,b,c) \mapsto (a,b,3ab-c).$$
This process is called a \emph{mutation}.
All solutions are obtained from $(1,1,1)$ by a sequence of mutations. 
The solutions lie at the vertices of an infinite tree, where two vertices are joined by an edge if they are related by a single mutation.
Each vertex has degree $3$, and there is a natural action of $S_3$ on the tree obtained by permuting the variables.
The first few solutions are $(1,1,1)$, $(1,1,2)$, $(1,2,5)$, $(1,5,13)$, $(2,5,29)$.
Hence Theorem~\ref{thm-manetti} provides a very explicit description of all Manetti surfaces.

\subsection{Log canonical surfaces}

We show that the log canonical degenerations of the plane are precisely the Manetti surfaces and the elliptic cones of
degree~9. 

\begin{lem} \cite{Ma} \label{lem-exclocus}  Let $X$ be a normal rational surface which smoothes to $\bP^2$.
Let $\pi \colon \tilde{X} \map X$ be the minimal resolution of $X$.
Then, assuming $X$ is not isomorphic to $\bP^2$, there exists a birational morphism $\tilde{X} \map \bF_w$; fix one such morphism
$\mu$ with $w$ maximal. Let $p \colon \tilde{X} \map \bP^1$ denote the birational ruling induced by $\mu$, let 
$B$ be the negative section of $\bF_w$ and $B'$ its strict transform on $X$. Then the exceptional locus of $\pi$ consists of the 
curve $B'$ together with the components of the fibres of $p$ with self intersection at most $-2$, i.e.,
$$\Exc(\pi)=B' \cup \{ \Gamma \subset \tilde{X} \mbox{ $|$ } p_{\star}\Gamma=0 \mbox{ and } \Gamma^2 \le -2 \}.$$
Moreover, every degenerate fibre of $p$ contains a unique $(-1)$-curve.
\end{lem}

\begin{thm} \label{thm-ell_cone}
Let $X$ be a normal log canonical surface with $-K_X$ ample which admits a smoothing to $\bP^2$.
Then either $X$ is log terminal or $X$ is an elliptic cone of degree $9$.
\end{thm}
\begin{proof}
The smoothing of $X$ is necessarily $\bQ$-Gorenstein by Propostion~\ref{prop-rho},
so $K_X^2=K_{\bP^2}^2=9$. If $X$ is not rational, then $X$ is an elliptic cone by Theorem~\ref{thm-rationality},
and $X$ has degree $9$ since $K_X^2=9$.

Now suppose that $X$ is rational. 
We assume that $X$ has a strictly log canonical singularity and obtain a contradiction.
We first describe the rational strictly log canonical surface singularities (\cite{KM}, p. 112--116). 
The exceptional locus $E$ of the minimal resolution is a collection of smooth rational curves in one of the following configurations:
\begin{enumerate}
\item[($\alpha$)] $E=G_1 \cup G_2 \cup F_1 \cup \cdots \cup F_k \cup G_3 \cup G_4$, where $F_1 \cup \cdots \cup F_k$ is a chain of smooth rational curves and 
$G_1,\cdots,G_4$ are $(-2)$-curves. The curves $G_1$ and $G_2$ each intersect $F_1$ in a single node and similiarly $G_3$ and $G_4$ each intersect $F_k$.
\item[($\beta$)] $E=F \cup G^1 \cup G^2 \cup G^3$, where $F$ is a smooth rational curve, each $G^i$ is a chain of smooth rational curves
$G^i_1 \cup \cdots \cup G^i_{k(i)}$ and the end component $G^i_1$ intersects $F$ in a single node.
\end{enumerate}

We use the notation and result of Lemma~\ref{lem-exclocus}.
Note first that $\mu \colon \tilde{X} \map \bF_w$ is an isomorphism over the negative section $B$ of $\bF_w$ by the maximality of $w$.  
Consider the MMP yielding $\mu \colon \tilde{X} \map \bF_w$ in a neigbourhood of a given degenerate 
fibre $f$ of the ruling $p \colon \tilde{X} \map \bP^1$.
At each stage we contract a $(-1)$-curve which meets at most 2 components of the fibre and is disjoint from $B'$.
By the Lemma, we have
$$\Exc(\pi)=B' \cup \{ \Gamma \subset \tilde{X} \mbox{ $|$ } p_{\star}\Gamma=0 \mbox{ and } \Gamma^2 \le -2 \}.$$
This set decomposes into the exceptional loci of the minimal resolutions of the log canonical singularity and some singularities of type
$\frac{1}{n^2}(1,na-1)$. Let $E$ denote the connected component of the exceptional locus contracting to the log canonical singularity.
Then $E$ has a component $C$ meeting 3 other components of $E$ --- we call such a curve a \emph{fork} of $E$.
We can now describe the form of the degenerate fibre $f$.  
Suppose first that $f$ contains a fork of $E$. Then we have a decompostion $f=P \cup \Gamma \cup Q \cup C \cup R \cup S,$ where
\begin{enumerate}
\item The curve $\Gamma$ is the unique $(-1)$-curve contained in $f$ and $C$ is a fork of $E$.
\item The curve $P$ is either empty or a chain of smooth rational curves, with one end component meeting $\Gamma$, 
which contracts to a singularity of type $\frac{1}{n^2}(1,na-1)$.
\item The curves $Q$, $R$ and $S$ are non-empty configurations of smooth rational curves such that $Q$ connects $\Gamma$ to $C$ 
and $S$ connects $C$ to $B'$ while $R$ intersects only $C$.
\end{enumerate}
Let $A'$ denote the component of $f$ meeting $B'$, i.e., the strict transform of the corresponding fibre of $\bF_w$.
Note that $A'$ cannot be a fork of $E$ since $f$ contains only one $(-1)$-curve; in particular $S$ is non-empty as claimed.
Also, in the MMP $\tilde{X} \map \cdots \map \bF_w$, we contract the components of $f-A'$ in the following order: 
$\Gamma$, $P \cup Q$, $C$, $R \cup S -A'$.

Suppose now that $f$ does not contain a fork of $E$. Then we have a decompositon $f=P \cup \Gamma \cup Q$, where
\begin{enumerate}
\item The curve $\Gamma$ is the unique $(-1)$-curve contained in $f$.
\item The curve $P$ is either empty or a chain of smooth rational curves, with one end component meeting $\Gamma$, 
which contracts to a singularity of type $\frac{1}{n^2}(1,na-1)$.
\item The curve $Q$ is a non-empty chain of curves meeting $\Gamma$, with one end component meeting $B'$.
\end{enumerate}

Suppose that $E$ has the form ($\alpha$) above.
If one of the forks $F_1, F_k$ of $E$ is contained in a degenerate fibre $f$, then we can write $f=P \cup \Gamma \cup Q \cup C \cup R \cup S$ as above where, without loss of generality, $Q=G_1$, $C=F_1$, $R=G_2$ and 
$S=F_2 \cup \cdots \cup F_l$ for some $l<k$.
Note that $E$ cannot contain the remaining fork $F_k$ of $E$ since then $F_k=A'$, a contradiction. 
Considering the MMP $\tilde{X} \map \cdots \map \bF_w$ again, we deduce that the curves in the chain $P$ have self-intersections $-3,-2,\ldots,-2$.
Thus $P$ contracts to an $\frac{1}{2r+1}(1,r)$ singularity, where $r$ is the length of the chain. But this is not of type $\frac{1}{n^2}(1,na-1)$,
a contradiction. Hence $P$ is empty. It follows that the curves in the chain 
$S=F_2 \cup \cdots \cup F_l$ have self-intersections $-3,-2,\ldots,-2,-1$  
if $l>2$. But $F_l^2 \le -2$, hence $l=2$ and $F_2^2=-2$. 
On the other hand, if the fork $F_1$ is not contained in a degenerate fibre, then $F_1=B'$. 
Then there is a degenerate fibre $f$ of the second form $P \cup \Gamma \cup Q$ with $Q=G_1$, a $(-2)$-curve.
It follows that the chain $P$ is a single $(-2)$-curve, which contracts to a $\frac{1}{2}(1,1)$ singularity, a contradiction.
Combining our results, we deduce that $k=5$, there are two fibres of the form $\Gamma \cup G_1 \cup G_2 \cup F_1 \cup F_2$ as above and $F_3=B'$.
There are no further degenerate fibres. We compute that $w=11$ using $K_X^2=9$.

We claim that the surface $X$ constructed above does not admit a smoothing to $\bP^2$.
Let $Z \map X$ be a local canonical covering of the singularity.
Then $Z \map X$ is a $\mu_2$ quotient and $Z$ has a cusp singularity.
The minimal resolution of $Z$ has exceptional locus a cycle of smooth rational curves with self-intersections $-2,-2,-2,-11,-2,-2,-2,-11$.
Suppose there exists a smoothing of $X$ to $\bP^2$, then we obtain a smoothing of $Z$ by taking the canonical covering.
Let $M$ denote the Milnor fibre of the smoothing of $Z$ and let $\mu_-$ denote the number of negative entries in a diagonalisation 
of the intersection form on $H^2(M,\bR)$. Then
$$\mu_- = 10 h^1(\cO_{\tilde{Z}}) + K_{\tilde{Z}}^2+b_2(\tilde{Z})-b_1(\tilde{Z})$$
where $\tilde{Z}$ is the minimal resolution of $Z$ \cite{St2}. In our case we calculate $\mu_- = 10-18+8-1=-1$, a contradiction.

Suppose now that $E$ has the form $(\beta)$. We first describe $E$ in more detail.
The chains $G^1$, $G^2$ and $G^3$ can be contracted to yield a partial resolution $\phi \colon \hat{X} \map X$; 
write $\hat{F}$ for the image of $F$ on $\hat{X}$.
The chains $G^i$ contract to singularities of the pair $(\hat{X},\hat{F})$ of type $(\frac{1}{r}(1,a),\Delta)$.
Let $r_1$, $r_2$ and $r_3$ be the indices of these singularities, then $\sum \frac{1}{r_i}=1$.
For $X$ is assumed to be strictly log canonical, hence $K_{\hat{X}}=\phi^{\star}K_X-\hat{F}$ or, equivalently,
$$0=(K_{\hat{X}}+\hat{F})\hat{F} = -2 + \sum (1-\frac{1}{r_i})=1-\sum\frac {1}{r_i}.$$
So $(r_1,r_2,r_3)=(2,3,6)$, $(2,4,4)$ or $(3,3,3)$. In particular, each chain $G^i$ is either a single smooth rational curve of self-intersection $-r_i$
or a chain of $(-2)$-curves of length $(r_i-1)$. 

We claim that the fork $F$ of $E$ cannot be contained in a fibre $f$.
It is enough to show that $w$ is greater than $2$.
For ${B'}^2=B^2=-w$, so if $w>2$ then $F=A'$ by the description of $E$ above, a contradiction.
Define an effective $\bQ$-divisor $\tilde{C}$ on $\tilde{X}$ by $K_{\tilde{X}}+\tilde{C}=\pi^{\star}K_X$ and let $C_1$ be the image of $\tilde{C}$
on $\bF_w$. Then
$$(K_{\bF_w}+C_1)^2 > (K_{\tilde{X}}+\tilde{C})^2=K_X^2=9.$$
Since $\mu$ is an isomorphism over the negative section $B$ of $\bF_w$, we have
$$(K_{\bF_w}+C_1)B=(K_{\tilde{X}}+\tilde{C})B'=\pi^{\star}K_X \cdot  B'=0.$$ 
So $K_{\bF_w}+C_1\sim \lambda (B+wA)$, where $A$ is a fibre of $\bF_w/\bP^1$.
Here $\lambda=-2+m$ where $m$ is the multiplicity of $B'$ in $\tilde{C}$ and $0 \le m \le 1$ since $X$ is log canonical. Hence
$$9 < (K_{\bF_w}+C_1)^2 =\lambda^2(B+wA)^2=\lambda^2w \le 4w,$$ 
so $w$ is greater than $2$ as required.

Thus $F=B'$ and there are 3 degenerate fibres $f_1$, $f_2$, and $f_3$ of the second form $P \cup \Gamma \cup Q$, where $Q=G^1$, $G^2$ and $G^3$ 
respectively. 
Recall that $G^i$ is either a single smooth curve of self-intersection $-r_i$ or a chain of $(-2)$-curves of length $(r_i-1)$.
If the fibre $f_i$ is a chain, we deduce that $P$ is either a chain of $(-2)$-curves of length $(r_i-1)$ 
or a single smooth curve of self-intersection $-r_i$ respectively. 
Since $P$ contracts to a singularity of type $\frac{1}{n^2}(1,na-1)$, it follows that $r_i=4$ and $P$ is a $(-4)$-curve.
If $f_i$ is not a chain, we find that $Q$ is a chain of three $(-2)$-curves, the $(-1)$-curve $\Gamma$ meets the middle component of $Q$ and $P$ is empty,
hence again $r_i=4$. So $r_i=4$ for each $i$, contradicting the description of $E$ above.
\end{proof} 

\begin{rem}
Conversely, it is well-known that an elliptic cone of degree $9$ admits a smoothing to $\bP^2$ \cite{Pi}.
\end{rem}

\section{The type B surfaces} \label{typeB}

We give necessary and sufficient conditions for a surface of type B to admit a smoothing to the plane.
Together with the description of the Manetti surfaces in Section~\ref{normal}, this completes the finer classification of
the surfaces appearing in stable pairs of degree $d$ not a multiple of $3$.

\begin{thm}\label{thm-B}
Let $X$ be a surface of type $B$. Then $X$ admits a $\bQ$-Gorenstein smoothing to $\bP^2$ if and only if the following conditions are 
satisfied:
\begin{enumerate}
\item The surface $X$ has singularities of the following types:
\begin{enumerate}
\item $\frac{1}{n^2}(1,na-1)$ where $(a,n)=1$.
\item $(xy=0) \subset \frac{1}{r}(1,-1,a)$ where $(a,r)=1$.
\end{enumerate}
Moreover there are at most $2$ singularities of type (b) of index $r$ greater than $1$.
\item $K_X^2=9$
\item Either (i) $\rho(X_1)=\rho(X_2)=1$ or (ii) $\{\rho(X_1),\rho(X_2)\}=\{1,2\}$.
\end{enumerate}
Moreover, in this case, $X$ has unobstructed $\bQ$-Gorenstein deformations.
\end{thm}

\begin{proof}
We first prove that the conditions are necessary.
The surface $X$ has only singularities of types (a) and (b) by Propositions~\ref{prop-smoothableslt} and 
\ref{prop-singsP2}. There are at most two singularities of type (b) by Lemma~\ref{lem-divcontr}.
We have $K_X^2=9$ since $X$ admits a $\bQ$-Gorenstein smoothing to $\bP^2$. Finally the Picard numbers of the components of $X$ are
as described in (3) by Proposition~\ref{prop-rho}.

Now suppose that $X$ satisfies the conditions above. We use the $\bQ$-Gorenstein deformation theory developed in Section~\ref{QG} to 
prove the existence of a smoothing. We first show that the $\bQ$-Gorenstein deformations of $X$ are unobstructed. 
It is enough to show that $H^p(\cT_{QG,X}^q)=0$ for $p+q=2$.
A local canonical covering of $X$ is a hypersurface, hence $H^0(\cT^2_{QG,X})=0$. 
The sheaf $\cT^1_{QG,X}$ is supported on the singular locus of $X$, which consists of the double curve
$\Delta$ together with some isolated points. The curve $\Delta$ is smooth and rational; let
$i \colon \bP^1 \inj X$ denote the inclusion of the $\Delta$. Near $\Delta$ the sheaf 
$\cT^1_{QG,X}$ equals either $i_{\star}\cO_{\bP^1}(1)$ or $i_{\star}\cO_{\bP^1}$ by Lemma~\ref{lem-T^1} and Lemma~\ref{lem-Delta^2},
where the two cases correspond to cases (i) and (ii) of condition (3) respectively. Hence in particular $H^1(\cT^1_{QG,X})=0$.
Finally $H^2(\cT^0_{QG,X})=H^2(\cT_X)=0$ by Lemma~\ref{lem-H^2(T_X)}.

We now construct a smoothing of $X$; we first construct an appropriate first order deformation of $X$ and then extend 
it to obtain a smoothing. 
Let $P \in X$ be a point of type $\frac{1}{n^2}(1,na-1)$, then
$$P \in X \cong (xy-z^n=0) \subset \frac{1}{n}(1,-1,a).$$  
Locally at $P$, the sheaf $\cT^1_{QG,X}$ equals the skyscraper sheaf $k(P)$ and a non-zero section
corresponds to a first order deformation of the form
$$(xy-z^n+t=0) \subset \frac{1}{n}(1,-1,a) \times \Spec (k[t]/(t^2)).$$
At the double curve $\Delta \cong \bP^1 \stackrel{i}{\inj} X$, 
the sheaf $\cT^1_{QG,X}$ equals either $i_{\star}\cO_{\bP^1}$ or $i_{\star}\cO_{\bP^1}(1)$.
Hence we may pick a section $s$ of $\cT^1_{QG,X}$ which is either nowhere zero on $\Delta$ or has a unique zero at $Q \in \Delta$, 
where $Q \in X$ is a normal crossing point. The section $s$ corresponds to a first order deformation of a neighbourhood of $\Delta$ in 
$X$ which is locally of the form 
$$(xy+t=0) \subset \frac{1}{r}(1,-1,a) \times \Spec (k[t]/(t^2))$$
away from the zeroes of $s$ and of the form
$$(xy+zt=0) \subset \bA^3 \times \Spec (k[t]/(t^2))$$
at a zero. Since $H^2(\cT_{QG,X})=0$, we can lift a section $s \in H^0(\cT^1_{QG,X})$ to an element of $T^1_{QG,X}$, so there is a 
global first order infinitesimal deformation of $X$ which is locally of the forms described above.
Given such a first order deformation of $X$, we can extend it to a $\bQ$-Gorenstein 
deformation $\cX/T$ over the germ of a curve since $\bQ$-Gorenstein deformations of $X$ are unobstructed.
Then the general fibre of $\cX/T$ is a smooth del Pezzo surface such that $K^2=9$, hence is isomorphic to $\bP^2$.
\end{proof}

\begin{lem}\label{lem-T^1} (\cite{Has}, Proposition~3.6)
Let $X$ be a surface with two normal irreducible components meeting in a smooth curve $\Delta$.
Suppose that $X$ has only singularities of the form $(xy=0) \subset \frac{1}{r}(1,-1,a)$ at $\Delta$.
Then, in a neighbourhood of $\Delta$,
$$\cT^1_{QG,X} \cong \cO_{\Delta}(\Delta_1 |_{\Delta}+\Delta_2 |_{\Delta}).$$
Here $\Delta_i$ is the restriction of $\Delta$ to $X_i$ and we calculate $\Delta_i |_{\Delta}$ by moving $\Delta_i$ on $X_i$ and
restricting to $\Delta$; thus $\Delta_i |_{\Delta}$ is a $\bQ$-divisor on $\Delta$ which is well defined modulo linear equivalence.
The sum $\Delta_1 |_{\Delta}+ \Delta_2|_{\Delta}$ is a $\bZ$-divisor on $\Delta$. In particular, the sheaf 
$\cT^1_{QG,X}$ is a line bundle on $\Delta$ of degree $\Delta_1^2+\Delta_2^2$.
\end{lem}

\begin{lem}\label{lem-Delta^2}
Let $X$ be a surface of type $B$ satisfying the conditions of Theorem~\ref{thm-B}.
Then
\begin{eqnarray*}
\Delta_1^2 +\Delta_2^2=\left\{  \begin{array}{cl}
                                        1 & \mbox{ if } \rho(X_1)=\rho(X_2)=1\\
                                        0 & \mbox{ if } \{\rho(X_1),\rho(X_2)\}=\{1,2\}
                                \end{array} \right.
\end{eqnarray*} 
\end{lem}
\begin{proof}
Let $\tilde{X_i} \map X_i$ be the minimal resolution of the component $X_i$ of $X$ for $i=1$ and $2$.
Then 
$$K_{\tilde{X_i}}^2+\rho(\tilde{X_i})=10$$
for each $i$ by Noether's formula and 
$$K_{\tilde{X_1}}^2+K_{\tilde{X_2}}^2+\rho(\tilde{X_1})+\rho(\tilde{X_2})=
K_{X_1}^2+K_{X_2}^2+\rho(X_1)+\rho(X_2)+ 4 \sum (1-\frac{1}{r_j}),$$
where the $r_j$ are the indices of the non-Gorenstein singularities of $X$ at $\Delta$ 
(cf. Theorem~\ref{thm-B*}). Thus
$$K_{X_1}^2+K_{X_2}^2 = 20 -(\rho(X_1)+\rho(X_2))-4 \sum (1-\frac{1}{r_j}).$$
The condition $K_X^2=9$ may be rewritten 
$$(K_{X_1}+\Delta_1)^2+(K_{X_2}+\Delta_2)^2 = 9.$$
Finally,
$$K_{X_i}\Delta_i+\Delta_i^2=-2+ \sum (1-\frac{1}{r_j})$$
for each $i$ by adjunction. Solving for $\Delta_1^2+\Delta_2^2$ we obtain $\Delta_1^2+\Delta_2^2= 3 - (\rho(X_1)+\rho(X_2))$.
\end{proof}

\begin{lem}\label{lem-H^2(T_X)}
Suppose $X$ is a surface of type B which satisfies the conditions of Theorem~\ref{thm-B}. Then $H^2(\cT_X)=0$.
\end{lem}

\begin{proof}
We have an exact sequence
$$ 0 \map \cO_{X_1}(-\Delta_1) \oplus \cO_{X_2}(-\Delta_2) \map \cO_X \map \cO_{\Delta} \map 0.$$
Applying the functor $\cHom_{\cO_X}(\Omega_X, \cdot)$, we obtain the exact sequence
$$ 0 \map \cT_{X_1}(-\Delta_1) \oplus \cT_{X_2}(-\Delta_2) \map \cT_X \map \cHom_{\cO_{\Delta}}(\Omega_X|_{\Delta}, \cO_{\Delta}).$$
Thus we have an inclusion $\cT_{X_1}(-\Delta_1) \oplus \cT_{X_2}(-\Delta_2) \inj \cT_X$ with cokernel supported on $\Delta$.
It follows that the map $H^2(\cT_{X_1}(-\Delta_1)) \oplus H^2(\cT_{X_2}(-\Delta_2)) \map H^2(\cT_X)$ is surjective.
So it is enough to show that $H^2(\cT_{X_i}(-\Delta_i))=0$ for $i=1$ and $2$.

Let $(Y,C)$ denote one of the pairs $(X_i,\Delta_i)$.
By Serre duality, 
$$H^2(\cT_{Y}(-C)) \cong \Hom(\cT_{Y}(-C),\cO_{Y}(K_{Y}))^{\vee}
= \Hom(\cT_{Y},\cO_{Y}(K_{Y}+C))^{\vee}.$$
We claim that $\cO_Y(-K_{Y}-C)$ has a nonzero global section. Assuming this,
$$\Hom(\cT_{Y},\cO_{Y}(K_{Y}+C)) \inj \Hom(\cT_{Y},\cO_{Y})=H^0(\Omega_{Y}^{\vee\vee}).$$
Now, letting $\pi : \tilde{Y} \map Y$ be the minimal resolution, we have 
$\Omega_{Y}^{\vee\vee}=\pi_{\star}\Omega_{\tilde{Y}}$ since $Y$ has only quotient singularities 
(\cite{St1}, Lemma~1.11). Thus $h^0(\Omega_{Y}^{\vee\vee})=h^0(\Omega_{\tilde{Y}})=h^1(\cO_{\tilde{Y}})=0$. 
So $H^2(\cT_{Y}(-C))=0$ as required.

It remains to show that $\cO_Y(-K_{Y}-C)$ has a nonzero global section.
Consider the exact sequence
$$0 \map \cO_Y(-K_Y-2C) \map \cO_Y(-K_Y-C) \map \cO_C(-K_Y-C) \map 0.$$
Now $h^1(\cO_Y(-K_Y-2C))=h^1(\cO_Y(2K_Y+2C))=0$ by Serre duality and Kodaira vanishing (recall that $Y$ is log terminal
and $-(K_Y+C)$ is ample). So it is enough to show that $\cO_C(-K_Y-C)$ has a nonzero global section.
A local calculation shows that $\cO_C(-K_Y-C) \cong \cO_C(-K_C-S)$, where $S$ is the sum of the singular points of $Y$ 
lying on $C$. Now $C$ is isomorphic to $\bP^1$, and there are at most $2$ singular points of $Y$ on $C$ by assumption, 
thus $\deg(-K_C-S) \ge 0$ and $\cO_C(-K_Y-C)$ has a nonzero global section as required.
\end{proof}

\begin{lem}\label{lem-divcontr}
Let $X$ be a surface of type B that admits a $\bQ$-Gorenstein smoothing to $\bP^2$. 
Then $X$ has at most two singularities of the form
$$(xy=0) \subset \frac{1}{r}(1,-1,a)$$
where the index $r$ is greater than $1$.
\end{lem}
\begin{proof}
Suppose that $X$ is a counter-example and let $\cX/T$ be a $\bQ$-Gorenstein smoothing of $X$ to $\bP^2$.
First assume that $\rho(X_1)=1$ and $\rho(X_2)=2$. We claim that there is a Mori contraction $f \colon \cX \map \cY/T$
with exceptional locus $X_1$. Assuming this, we deduce that the special fibre $Y$ of $\cY/T$ has a log terminal singularity 
such that the exceptional locus of the minimal resolution has a `fork', i.e., 
there exists an exceptional curve meeting 3 other exceptional curves. But, by 
Proposition~\ref{prop-singsP2}, the only possible log terminal singularities on $Y$ are cyclic quotient singularities, so the 
exceptional locus of the minimal resolution is a chain of curves, a contradiction.

We now prove the existence of the contraction $f$. 
We have $\Delta^2_1+\Delta^2_2=0$ by Lemma~\ref{lem-Delta^2} and $\rho(X_1)=1$ by assumption, hence  $\Delta_1^2>0$ and $\Delta_2^2<0$.
Thus  $\Delta_2$ generates an extremal ray on $X_2$. It follows that $\Delta$ generates an extremal ray on $\cX/T$.
The divisor $-K_{\cX}$ is relatively ample, so in particular $K_{\cX} \Delta < 0$ and there is a corresponding contraction 
$f \colon \cX \map \cY/T$.
The exceptional locus of $f$ is the divisor $X_1$ since $\Delta_1$ generates the group $N_1(X_1)$ of 1-cycles on $X_1$.

Similiarly, if $\rho(X_1)=\rho(X_2)=1$, then $\cX$ is not $\bQ$-factorial and there is a $\bQ$-factorialisation 
$\alpha \colon \hat{\cX} \map \cX/T$,
where the special fibre $\hat{X}$ of $\hat{\cX}$ has components $\hat{X}_1 \cong X_1$ and $\hat{X}_2$, a blowup of $X_2$.
Then there is a Mori contraction $f \colon \hat{\cX} \map \cY/T$ with exceptional locus $\hat{X_1}$; 
we obtain a contradiction as above. We construct the $\bQ$-factorialisation $\alpha$ explicitly below.
Let $P \in X$ be a point at which $\cX$ is not $\bQ$-factorial, then necessarily $P \in \Delta$ and, working locally analytically at 
$P \in \cX/T$, the family $\cX/T$ is of the form
$$(xy+t^kg(z^r,t)=0) \subset \frac{1}{r}(1,-1,a,0),$$
where $t$ is a local parameter at $0 \in T$ and $g(z^r,t) \in m_{\cX,P}$, $t \nd g(z^r,t)$. Let $X_1=(x=t=0)$ and $X_2=(y=t=0)$.
If $r=1$, let $\alpha \colon \hat{\cX} \map \cX$ be the blowup of $(x=g=0) \subset \cX$.
Then, writing $u=g/x$ and $v=x/g$, the 3-fold $\hat{X}$ has the following affine pieces:
$$(vy+t^k=0) \subset \bA^4_{v,y,z,t}$$
$$(xu=g(z^r,t)) \subset \bA^4_{x,u,z,t}$$
Thus $\hat{X}_1$ is isomorphic to $X_1$ and the morphism $\hat{X}_2 \map X_2$ contracts a smooth rational curve to the point 
$P \in X_2$.
If $r>1$, we obtain $\alpha$ as the quotient of the above construction applied to the canonical covering of $\cX$.
Finally $\hat{\cX}$ is $\bQ$-factorial since $\rho(\hat{X}_1)=1$ and $\rho(\hat{X}_2)=2$, cf. Proposition~\ref{prop-rho}.
\end{proof}

\section{The singularities of $D$ and the relation to GIT} \label{GIT}

If $(X,D)$ is a stable pair then the pair $(X,(\frac{3}{d}+\epsilon)D)$ is slc for $0 < \epsilon \ll 1$.
We show that, in the case $X=\bP^2$, this condition is a natural strengthening of the GIT stability condition.
Roughly speaking, it is the weakest local analytic condition on $D$ which contains the GIT stability condition.
This statement is made precise in Propositions~\ref{prop-stability} and \ref{prop-GIT}.

\begin{defn}
Let $P \in X$ be the germ of a smooth surface and $D$ a divisor on $X$.
Suppose given a choice of coordinates $x$, $y$ at $P \in X$ and weights $(m,n) \in \bN^2$.
Write $D=(f(x,y)=0)$ and $f(x,y)=\sum a_{ij}x^iy^j$.
The \emph{weight} $\wt(D)$ of $D$ is given by
$$\wt(D)=\min \{ mi+nj \mbox{ $|$ } a_{ij} \neq 0 \}.$$
\end{defn}

\begin{prop} \label{prop-stability}
Let $D$ be a plane curve of degree $d$.
Then $(\bP^2,D)$ is a stable pair if and only if for every point $P \in \bP^2$, choice of analytic
coordinates $x$, $y$ at $P$ and weights $(m,n)$, we have $$\wt(D) < \frac{d}{3}(m+n).$$
\end{prop}
\begin{proof}
Given a smooth surface $X$ and $B$ a $\bQ$-divisor on $X$, to verify that $(X,B)$ is log canonical 
it is sufficient to check that, for each weighted blowup 
$$f \colon E \subset Y \map P \in X$$ 
of a point $P \in X$, we have $a(E,X,B) \ge -1$. Here $a=a(E,X,B)$ is the discrepancy defined by the equation
$$K_Y+B'=f^{\star}(K_X+B)+aE.$$ 
Putting $X=\bP^2$ and $B=(\frac{3}{d}+\epsilon)D$ yields the criterion above.
\end{proof}

\begin{defn}
We say that coordinates $x$, $y$ at a point $P \in \bP^2$ are \emph{linear} if there is a choice of 
homogeneous coordinates $X_0$, $X_1$, $X_2$ on $\bP^2$ such that $x=X_1/X_0$ and $y=X_2/X_0$.
\end{defn}

\begin{prop} \label{prop-GIT}
Let $D$ be a plane curve of degree $d$.
Then $D \inj \bP^2$ is GIT stable if and only if for every point $P \in \bP^2$, choice of linear coordinates 
$x$, $y$ at $P$ and weights $(m,n)$, we have
$$\wt(D) < \frac{d}{3}(m+n).$$
\end{prop}
\begin{proof}
This is the usual numerical criterion for GIT stability \cite{Mu}, restated in a form analogous to 
Proposition~\ref{prop-stability}. 
\end{proof}

\begin{ex}
We give an example of a curve $D \subset \bP^2$ such that $D$ is GIT stable
but $(\bP^2,D)$ is not a stable pair. The curve $D$ is a 
quintic curve with a singularity $P \in D$ of type 
$(y^2+x^{13}=0) \subset \bC^2$. 
To prove the existence of such a curve, pick analytic coordinates $x$, $y$ at
$P=(1:0:0) \in \bP^2$. Let $F$ be a homogenenous polynomial of degree 5, and
write $F/X_0^5$ as a power series $f(x,y)$ in $x$ and $y$.
Quintics depend on 20 parameters, hence we may choose $F$ so that 
the coefficients of
$1, x,\ldots, x^{12}, y, xy,\ldots, x^6y$ in $f$ vanish.
Then, for sufficiently generic choice of $x$ and $y$,
$f(x,y)=\alpha y^2 +\beta x^{13} + \cdots$, where $\alpha \neq 0$, $\beta \neq 0$
and $\cdots$ denotes terms of higher weight with respect to the weights 
$(2,13)$ of $x$ and $y$. In this case, the quintic curve $D=(F=0)$ has a 
singularity of the desired type at $P$.
Then $(\bP^2,D)$ is not a stable pair by Proposition~\ref{prop-stability} --- 
with respect to the weighting $(2,13)$ of $x$ and $y$ we have
$\wt(f) =26 > \frac{5}{3}(2+13)=25$. On the other hand, 
let $D' \map D$ be the resolution of 
the singularity $P \in D$ induced by a $(2,13)$ weighted blowup of 
$\bP^2$.
We compute that $p_a(D')=0$, hence $D'$ is a smooth rational curve and 
$D$ has no additional singular points. Thus $D$ is
GIT stable by (\cite{Mu}, p. 80).
\end{ex}

\section{Examples} \label{Examples}

We give the classification of stable pairs of degrees 4 and 5.
 
\begin{notn}
Given an embedding of a surface $Y$ in a weighted projective space $\bP$, we write $kH$ for a general curve in the linear system 
$|\cO_{Y}(k)|$. For each surface of type B, we use this notation to describe the inverse image of the double curve on each component.
\end{notn}

When we list the singularities of the surfaces $X$ 
we do not mention the normal crossing singularities $(xy=0) \subset \bA^3$.
Similiarly, when we list the possible singularities of $(X,D)$, we do not include the cases
where $X$ is smooth or normal crossing and the divisor $D$ is normal crossing.   

\subsection{Degree 4} Surfaces $X$:
\begin{eqnarray*}
\renewcommand{\arraystretch}{1.5}
\begin{array}{|l|l|l|} \hline
\mbox{Surface}          & \mbox{Double curve}   & \mbox{Singularities} \\ \hline
\bP^2                           &                       &                       \\
\bP(1,1,4)                      &                       & \frac{1}{4}(1,1)      \\
\bP(1,1,2) \cup \bP(1,1,2)      & H,H                   & (xy=0) \subset \frac{1}{2}(1,1,1) \\ \hline
\end{array}
\end{eqnarray*}
Allowed singularities of $(X,D)$:
\begin{eqnarray*}
\renewcommand{\arraystretch}{1.5}
\begin{array}{|l|l|} \hline
X       & D     \\ \hline
\bA^2_{x,y}                             & (y^2+x^3=0)\\
\frac{1}{4}(1,1)                &   0 \\
(xy=0) \subset \frac{1}{2}(1,1,1) & 0\\ \hline
\end{array}
\end{eqnarray*}

\subsection{Degree 5} 
Surfaces $X$:
\begin{eqnarray*}
\renewcommand{\arraystretch}{1.5}
\begin{array}{|l|l|l|} \hline
\mbox{Surface}                  & \mbox{Double curve}   & \mbox{Singularities} \\ \hline
\bP^2                                   &                       &                       \\
\bP(1,1,4)                              &                       & \frac{1}{4}(1,1)      \\
X_{26} \subset \bP(1,2,13,25)           &                       & \frac{1}{25}(1,4)     \\
\bP(1,4,25)                             &                       & \frac{1}{4}(1,1), \frac{1}{25}(1,4) \\
\bP(1,1,2) \cup \bP(1,1,2)              & H,H                   & (xy=0) \subset \frac{1}{2}(1,1,1) \\
\bP(1,1,5) \cup (X_6 \subset \bP(1,2,3,5))& H, 2H               & (xy=0) \subset \frac{1}{5}(1,-1,1) \\
\bP(1,1,5) \cup \bP(1,4,5)              & H, 4H         & \frac{1}{4}(1,1), (xy=0) \subset \frac{1}{5}(1,-1,1) \\ \hline
\end{array}
\end{eqnarray*}
Allowed singularities of $(X,D)$:
\begin{eqnarray*}
\renewcommand{\arraystretch}{1.5}
\begin{array}{|l|l|} \hline
X               & D     \\ \hline
\bA^2_{x,y}     & (y^2+x^n=0) \mbox{ for } 3 \le n \le 9\\
\bA^2_{x,y}     & (x(y^2+x^n)=0) \mbox{ for } n=2,3 \\
\frac{1}{4}(1,1) & (y^2+x^n=0) \mbox{ for } n=2,6 \\
(xy=0) \subset \frac{1}{2}(1,1,1) & (z=0) \\ 
\frac{1}{25}(1,4) & 0 \\
(xy=0) \subset \frac{1}{5}(1,-1,1) & 0 \\ \hline
\end{array}
\end{eqnarray*}

\vspace{0.5cm}

Note that $X_{26} \subset \bP(1,2,13,25)$ is the surface 
obtained from $\bP(1,4,25)$ by smoothing the $\frac{1}{4}(1,1)$ singularity. 
The smoothing can be realised inside $\bP(1,2,13,25)$.
To see this, let $k[U,V,W]$ be the homogeneous coordinate ring of $\bP(1,4,25)$ and consider the 2nd Veronese subring 
$k[U,V,W]^{(2)}$. By picking generators for this ring we obtain the embedding
\begin{eqnarray*}
\begin{array}{ccc}
\bP(1,4,25) & \stackrel{\sim}{\longrightarrow} & (XT=Z^2) \subset \bP(1,2,13,25)\\
(U,V,W)     & \longmapsto & (X,Y,Z,T)=(U^2,V,UW,W^2)
\end{array}
\end{eqnarray*}
Then the smoothing of the $\frac{1}{4}(1,1)$ singularity is given by
$$(XT=Z^2+tT^{13}) \subset \bP(1,2,13,25) \times \bA^1_t.$$
Similiarly $X_6 \subset \bP(1,2,3,5)$ is the surface obtained from $\bP(1,4,5)$ by smoothing the $\frac{1}{4}(1,1)$ singularity.

\subsection{Sketch of proof}
We describe two different ways to establish the classification of stable pairs of degrees $4$ and $5$ given above.
We note immediately that all the surfaces $X$ occurring are either Manetti surfaces or type B surfaces by 
Theorem~\ref{thm-degnotmult3}.

\subsubsection{The geometric method}
We first classify \emph{semistable} pairs of degree $d$ using the classification of Manetti surfaces $X$ 
(Theorem~\ref{thm-manetti_defns}) and the bound on the index of the singularities (Theorem~\ref{thm-index}).
The possible surfaces $X$ for $d=4$ are $\bP^2$ and $\bP(1,1,4)$, whereas for $d=5$ we have
$\bP^2$, $\bP(1,1,4)$, $X_{26} \subset \bP(1,2,13,25)$ and $\bP(1,4,25)$.

We now deduce the classification of the stable pairs of degree $d$ using the following result:
\begin{prop}
Every stable pair $(X,D)$ of type B has a smoothing $(\cX,\sD)/T$
which is obtained from a smoothing $(\cY,\sD_{\cY})/T$ of a semistable pair by a
divisorial extraction, possibly followed by a flopping contraction. 
Moreover the divisorial extraction 
$f \colon (\hat{\cX},\hat{\sD}) \map (\cY,\sD_{\cY})/T$ is crepant in the following sense:
$K_{\hat{\cX}}+\frac{3}{d}\hat{\sD}=f^{\star}(K_{\cY}+\frac{3}{d}\sD_{\cY})$.
\end{prop}
\noindent This is a special case of the `stabilisation process' described in the proof of Theorem~\ref{thm-stable}, 
which produces a smoothing of a stable pair from a smoothing of a semistable pair. 
The proof uses the explicit construction in the proof of Lemma~\ref{lem-divcontr}.
Restricting to the special fibre, we see that the centre $P \in Y \subset \cY$ 
of the divisorial contraction is a strictly log canonical 
singularity of the  pair $(Y,\frac{3}{d}D_Y)$. If $d=4$ we deduce that $(P\in Y,D) \cong (\bA^2,(y^2+x^4=0))$.
For $d=5$ there are three possibilities for $(P \in Y, D)$, namely $(\bA^2,(y^2+x^{10}=0))$, $(\bA^2,(x(y^2+x^4)=0))$ and
$(\frac{1}{4}(1,1),(y^2+x^{10}=0))$. The required divisorial extractions $f \colon \hat{\cX} \map \cY$ 
are then determined by \cite{Hac2}.
The special fibre $\hat{X}$ of $\hat{\cX}$ is $Y'+E$ where 
$Y'$ is the strict transform of $Y$ and $E$ is the exceptional divisor of $f$.
The map $Y' \map Y$ is a weighted blowup with respect to some analytic coordinates
$x$, $y$ at $P \in Y$ as above. It is important to note that, for example in the case $Y=\bP^2$, these coordinates are not 
necessarily `linear coordinates' $X_1/X_0$, $X_2/X_0$ corresponding to homogeneous coordinates $X_0,X_1,X_2$ on $\bP^2$.
Hence the global structure of the rational surface $Y'$ is a little more complicated than one might expect.
Finally, if there is a curve $\Gamma$ on $Y' \subset \hat{\cX}$ such that $K_{\hat{\cX}}\Gamma=0$ then 
there is a flopping contraction $\alpha \colon \hat{\cX} \map \cX$ with exceptional locus $\Gamma$; otherwise $\cX=\hat{\cX}$. 
Thus either $X$ is obtained from $\hat{X}$ by contracting the curve $\Gamma \subset Y'$ or $X=\hat{X}$.

\subsubsection{The combinatorial method}
This approach is carried out carefully in \cite{Hac1}.
We set up the following notation: given a stable pair $(X,D)$, let $(Y,C)$ be a component of the pair $(X^{\nu},\Delta^{\nu})$, 
where $X^{\nu}$ is the normalisation of $X$ and $\Delta^{\nu}$ is the inverse image of the double curve of $X$.
Let $\pi \colon \tilde{Y} \map Y$ be the minimal resolution of $Y$ and define an effective $\bQ$-divisor $\tilde{C}$ 
by the equation
$$K_{\tilde{Y}}+\tilde{C}=\pi^{\star}(K_Y+C).$$
Assuming $Y$ is not isomorphic to $\bP^2$, there exists a birational morphism $\tilde{Y} \map \bF_w$; fix one such morphism
$\mu$ with $w$ maximal and let $p \colon \tilde{Y} \map \bP^1$ denote the induced ruling.

We first use the bound on the index of the singularities of $X$ (Theorem~\ref{thm-index}) to write down a list of possible 
singularities of the pair $(Y,C)$. We deduce the possible forms of the connected components of the divisor $\tilde{C}$.
We then analyse how these can embed into the surface $\tilde{Y}$ relative to the ruling $p$  
(cf. proof of Theorem~\ref{thm-ell_cone}). 
We deduce a list of candidates for the pairs $(\tilde{Y},\tilde{C})$ and hence for the pairs
$(Y,C)$. Finally we glue these components together to obtain the list of surfaces $X$.

\appendix
\section{The relative $S_2$ condition} \label{app_S_2}
 
\begin{defn} \label{S_2}
Let $\cX/S$ be a flat family of slc surfaces and $\cF$ a coherent sheaf on $\cX$.
We say $\cF$ is \emph{$S_2$ over $S$} if $\cF$ is flat over $S$ and the fibre $\cF_s = \cF \otimes k(s)$ 
satisfies Serre's $S_2$ condition for each $s \in S$.
We say $\cF$ is \emph{weakly $S_2$ over $S$} if, for each open subscheme $i \colon \cU \inj \cX$
whose complement has finite fibres, we have $i_{\star}i^{\star}\cF=\cF$. 
\end{defn}

\begin{rem} The relative $S_2$ condition is stable under base change, but this is \emph{not} true for the weak 
relative $S_2$ condition.
\end{rem}

\begin{lem}\label{weakly_S_2}
Let $\cX/S$ be a flat family of slc surfaces.
Let $\cF$ be a sheaf on $\cX$ which is $S_2$ over $S$.
Then $\cF$ is weakly $S_2$ over $S$. 
\end{lem}

\begin{ex} The sheaf $\cO_{\cX}$ is $S_2$ over $S$, hence $i_{\star}\cO_{\cU}=\cO_{\cX}$ for $i \colon \cU \inj \cX$
as in \ref{S_2}. Also, the sheaf $\omega_{\cX/S}$ is $S_2$ over $S$
(since $\omega_{\cX/S}$ is flat over $S$ and has fibres $\omega_{\cX_s}$ which are $S_2$ by \cite{KM}, Corollary~5.69),
so $i_{\star}\omega_{\cU/S}=\omega_{\cX/S}$.
\end{ex}

\begin{proof}
Let $i \colon \cU \inj \cX$ be an open subscheme as in \ref{S_2} and
let $\cZ$ denote the complement of $\cU$ with its reduced structure.
We work locally at a closed point $P \in \cZ$, let $P \mapsto s \in S$.
The sheaf $\cF_s$ is $S_2$ by assumption, so there is a regular sequence
$x_s,y_s \in m_{\cX_s,P}$ for $\cF_s$ at $P$. 
Now $\cZ_s \inj \cX_s$ is a closed subscheme with support $P$,
hence, replacing $x_s,y_s$ by powers $x_s^{\nu},y_s^{\nu}$ if necessary, 
we may assume that they lie in the ideal of $\cZ_s$. Note that $x_s,y_s$ is still a
regular sequence for $\cF_s$ by \cite{Mat}, Theorem~16.1.
Lift $x_s,y_s$ to elements $x,y$ of the ideal of $\cZ$, then $x,y$ is a regular
sequence for $\cF$ at $P$ (\cite{Mat}, p. 177, Corollary to Theorem~22.5).
Equivalently, we have an exact sequence
$$0 \map \cF \stackrel{(y,-x)}{\map} \cF \oplus \cF \stackrel{(x,y)}{\map} \cF.$$
Consider the natural map  $\cF \map i_{\star}i^{\star} \cF$, write $K$ for the kernel and $C$ for the cokernel.
Then $K$ and $C$ have support contained in the set $\cZ$, so any given element of $K$ or $C$ is annihilated by some 
power of the ideal $\cI_{\cZ}$ of $\cZ$.
So, if $K \neq 0$, there exists $0 \neq g \in K$ such that $\cI_{\cZ} g = 0$, then $xg=yg=0$, 
contradicting the exact sequence above. Similiarly if $C \neq 0$, there exists $g \in i_{\star}i^{\star}\cF \backslash \cF$ 
such that
$\cI_{\cZ}g \subset \cF$. Again using the exact sequence above, since $(yg,-xg) \mapsto 0$ we obtain $(yg,-xg)=(yg',-xg')$
for some $g' \in \cF$; it follows that $g=g'$, a contradiction. Thus $K=C=0$, so the map $\cF \map i_{\star}i^{\star} \cF$
is an isomorphism as claimed.
\end{proof}

\begin{lem} \label{S_2-basics}
Let $\cX/S$ be a flat family of slc surfaces.
\begin{enumerate}
\item If $\cF$ and $\cG$ are coherent sheaves on $\cX$ and $\cG$ is weakly $S_2$ over $S$, then $\cHom(\cF,\cG)$ is weakly 
$S_2$ over $S$.
\item If $0 \map \cF' \map \cF \map \cF'' \map 0$ is an exact sequence of coherent sheaves on $\cX$ and $\cF'$ and $\cF''$ are 
weakly $S_2$ over $S$, then $\cF$ is weakly $S_2$ over $S$.
\item Let $\cZ/S$ be a flat family of slc surfaces and $\pi \colon \cZ \map \cX$ a finite map over $S$. If 
$\cF$ is a sheaf on $\cZ$ which is weakly $S_2$ over $S$ then $\pi_{\star}\cF$ is weakly $S_2$ over $S$  
\item Let $g \colon T \inj S$ be a closed subscheme and $g_{\cX} \colon \cX_T \inj \cX$ the 
corresponding closed subscheme of $\cX$. If $\cF$ is a sheaf on $\cX_T$ which is weakly $S_2$ over $T$ then 
${g_{\cX}}_{\star}\cF$ is weakly $S_2$ over $S$.  
\end{enumerate}
\end{lem}
\begin{proof}
Let $i \colon \cU \inj \cX/S$ be an open subscheme whose complement has finite fibres.
For $\cF$ a sheaf on $\cX$, let $\alpha_{\cF}$ denote the natural map $\cF \map i_{\star}i^{\star}\cF$;
thus $\cF$ is weakly $S_2$ if and only if $\alpha_{\cF}$ is an isomorphism for each $\cU$.
To prove (1), observe that the map $\alpha_{\cHom(\cF,\cG)} \colon \theta \mapsto i_{\star}i^{\star}\theta$
has inverse $\psi \mapsto \alpha_{\cG}^{-1} \circ \psi \circ \alpha_{\cF}$. For (2), consider the diagram
\begin{eqnarray*}
\begin{array}{ccccccccc}
0   & \map & \cF' & \map & \cF & \map & \cF'' & \map & 0 \\
    &      &  \da \rlap{$\alpha_{\cF'}$} &      & \da \rlap{$\alpha_{\cF}$} &      &  \da \rlap{$\alpha_{\cF''}$}  &      & \\
0   & \map & i_{\star}i^{\star}\cF'& \map & i_{\star}i^{\star}\cF & \map & i_{\star}i^{\star}\cF'' & & \\ 
\end{array}
\end{eqnarray*}
The rows are exact, and $\alpha_{\cF'}$ and $\alpha_{\cF''}$ are isomorphisms by assumption, 
hence $i_{\star}i^{\star}\cF \map  i_{\star}i^{\star}\cF''$ is surjective and $\alpha_{\cF}$ is an isomorphism, as required.
Parts (3) and (4) follow immediately from the definition of the weak $S_2$ property.
\end{proof}

\noindent Department of Mathematics, University of Michigan,\\ 
Ann Arbor, MI 48109, USA\\
phacking@umich.edu

\end{document}